\documentclass[smallcondensed,runningheads,twosided]{svjour3}    
\usepackage{amssymb,amsmath,amsfonts,mathtools}                  
\usepackage{mathabx}                                             
\journalname{Submitted}                                          

\smartqed

\usepackage{mathabx}
\newcommand{\id}{\mathrm{id}}

\newcommand{\R}{\mathbb{R}}

\newcommand{\C}{\mathbb{C}}

\newcommand{\pa}{\partial}
\newcommand{\fJ}{\mathfrak{I}}
\newcommand{\la}{\langle}
\newcommand{\ra}{\rangle}
\newcommand{\Hb}{\mathbb{H}}

\begin{document}
\title{\textsf{Evans function and Fredholm determinants}}
\author{\textsf{Issa Karambal \and Simon J.A. Malham}}
\authorrunning{\textsf{Karambal \and Malham}}
\institute{Issa Karambal \and Simon J.A Malham \at       
Maxwell Institute for Mathematical Sciences\\            
and School of Mathematical and Computer Science\\        
Heriot--Watt University, Edinburgh EH14 4AS, UK\\        
Tel.: +44-131-4513200\\                                  
Fax: +44-131-4513249\\                                   
\email{ikarambal@gmail.com}\\                            
\email{simonmalham@gmail.com}}                           

\date{Received: 26th November 2014}
\voffset=10ex                    

\maketitle

\begin{abstract}
We explore the relationship between the Evans function, transmission
coefficient and Fredholm determinant for systems of first order linear differential
operators on the real line. The applications we have in mind include
linear stability problems associated with travelling wave solutions to
nonlinear partial differential equations, for example reaction-diffusion
or solitary wave equations. The Evans function and transmission coefficient, 
which are both finite determinants, are natural tools for both analytic 
and numerical determination of eigenvalues of such linear operators. However, 
inverting the eigenvalue problem by the free state operator generates a natural 
linear integral eigenvalue problem whose solvability is determined through 
the corresponding infinite Fredholm determinant. The relationship between
all three determinants has received a lot of recent attention. We focus
on the case when the underlying Fredholm operator is a trace class perturbation
of the identity. Our new results include: 
(i) clarification of the sense in which the Evans function and transmission
coefficient are equivalent; and (ii) proof of the equivalence of the 
transmission coefficient and Fredholm determinant, in particular 
in the case of distinct far fields.
\keywords{Fredholm determinant \and Evans function \and travelling waves}
\subclass{47B10 \and 47G10 \and 34B27 \and 34L05}
\end{abstract}

\newpage

\section{\textsf{Introduction}}
Our goal is to establish the connection between the
Evans function, transmission coefficient and Fredholm determinant
associated with linear $n$th order eigenvalue problems on $\R$ of the form
\begin{equation*}
(\pa-A_0-V)Y=O.
\end{equation*}
Here $\pa$ is the derivative operator $\pa Y=Y'$  
and $A_0\colon\R\times\C\mapsto\C^{n\times n}$ and $V\colon\R\mapsto\C^{n\times n}$
are bounded multiplicative operators. We suppose that $V$ represents a 
perturbative potential function that decays to zero in the far
field of the domain $\R$, while $A_0$ generates a background or free
state. It is constant in the far field though the limits are not 
necessarily the same. We suppose further that $A_0$ depends linearly
on a spectral parameter $\lambda\in\C$. Indeed, large classes of eigenvalue
problems can be couched in the form above. The problem is to determine
those values of $\lambda$, eigenvalues, for which suitable integrable 
solutions $Y\in\C^n$ exist to the equation above. 
The Evans function and transmission coefficient are standard tools in this
endeavour. Away from the essential spectrum, and suitably scaled, they are 
analytic functions of the spectral parameter $\lambda$ whose zeros coincide with
eigenvalues. The multiplicity of the zeros coincide with the algebraic
multiplicity of the eigenvalues. Modulo a non-zero scalar factor that renders
it domain independent, the Evans function is the determinant of the square matrix
whose left block is $Y^-$ and right block $Y^+$. The columns of these two matrices 
are solutions to the differential equation above that decay to zero exponentially fast
in the left and right far fields, respectively. The Evans function measures the 
``distance from intersection'' of the subspaces spanned by the columns of $Y^-$ and $Y^+$.
The transmission coefficient which is also a determinant,
measures the degree to which the solutions $Y^-$, that 
decay to zero in the left far field, are orthogonal to the subspace that
is orthogonal to the subspace of solutions that decays to zero in the right far field.
Unwrapping the two orthogonality conditions explains why the Evans function
and transmission coefficient are essentially equivalent. 
We assume away from the essential spectrum $(\pa-A_0)^{-1}$ exists.
Then our eigenvalue problem can be expressed in the form 
$\bigl(\id-(\pa-A_0)^{-1}V\bigr)Y=O$, or, with $V=U|V|$ representing the 
polar decomposition of $V$ and setting $\phi\coloneqq|V|^{1/2}Y$, in 
Birman--Schwinger form
\begin{equation*}
\bigl(\id-|V|^{1/2}(\pa-A_0)^{-1}U|V|^{1/2}\bigr)\phi=O.
\end{equation*}
From this perspective, we again seek values of the spectral parameter $\lambda\in\C$
for which solutions to this problem that decay to zero in the far field exist.
The natural underlying Hilbert space is $L^2(\R;\C^n)$. For the applications
we have in mind, establishing that $|V|^{1/2}(\pa-A_0)^{-1}U|V|^{1/2}$ 
is a Hilbert--Schmidt compact operator on this space is relatively straightforward. 
However herein we focus on the case when it 
is a trace class operator, i.e.\/ a nuclear operator. With this property,
zeros of the Fredholm determinant of $\id-|V|^{1/2}(\pa-A_0)^{-1}U|V|^{1/2}$
coincide with eigenvalues. Thus we come to the central issue. In what
sense are the Evans function, transmission coefficient and 
Fredholm determinant related? 
Let us briefly outline what has already been established.

The Evans function was first proposed by Evans~\cite{JE75}, while
Alexander, Gardner \& Jones~\cite{AGJ90} established it as a geometric
tool for stability analysis. Subsequently it has become a standard tool
in analytical and numerical studies of the stability of travelling waves;
see the review papers featuring the Evans function by Sandstede~\cite{BS02} 
and Kapitula~\cite{K}. The Evans function is also called the miss-distance
function, see Greenberg \& Marletta~\cite{GM}. It is also a generalization 
of the Wronskian and Jost function. The transmission coefficient has it origins 
much further back in the mathematical literature. Its connection to the Evans function, 
though trivial in the scalar case, can be found in Swinton~\cite{JS92} 
and Bridges \& Derks~\cite{BD99} for higher order problems. The Fredholm
determinant for determining the solvability of linear integral equations
was introduced by Fredholm~\cite{Fred}. Its connection to the transmission
coefficient goes back to Jost \& Pais~\cite{JP51}. Simon~\cite{BS04,BS05}
computes the explicit relationship between the Fredholm determinant and
Wronskian for some example scalar Schr\"odinger operators; also see
Kapitula \& Sandstede~\cite{KS04}. However more generally, 
Gesztesy \& Makarov~\cite{GK04} showed that for operators with
semi-separable kernels, their Fredholm and $2$-modified Fredholm determinants
can be reduced to the determinant of finite rank operators,
potentially useful for the evaluation of such Fredholm determinants.
Gesztesy, Latushkin \& Makarov~\cite{GL07} then established the
connection between the Evans function and a $2$-modified Fredholm determinant.
They also gave a coordinate free definition of the Evans function as
a ratio of the perturbed and unperturbed versions of the function.
The $2$-modified Fredholm determinant is relevant for their equivalence
results as systems of first order operators generate operators 
$|V|^{1/2}(\pa-A_0)^{-1}U|V|^{1/2}$ which are Hilbert--Schmidt class 
and in general not trace class. When such operators are trace class, 
the Fredholm determinant is the natural object in the equivalence result. 
Indeed, systems of Schr\"odinger operators represent
an explicit example, see Gesztesy, Latushkin \& Zumbrun~\cite[Section~4]{GYZ08}.

Our goal herein is to establish a unified picture of the relationship 
between the Evans function, transmission coefficient and Fredholm determinant.
We focus on those systems of first order operators for which 
$|V|^{1/2}(\pa-A_0)^{-1}U|V|^{1/2}$ is trace class and 
the matrix trace of the matrix perturbation potential $V$ is zero. 
By considering this subclass of first order operators we gain a degree of 
clarity and directness.
To begin with we assume $A_0$ is constant, but in our final 
main Section~\ref{sec:distinctfarfields} we assume distinct far field limits 
for $A_0$ which is therefore no longer constant. The free Evans function
and free transmission coefficients are the corresponding quantities associated
with the operator $\pa-A_0$. What we achieve in this paper is as follows, we:
\begin{enumerate}
\item Provide practical tests to determine 
when $|V|^{1/2}(\pa-A_0)^{-1}U|V|^{1/2}$ is trace class. These follow 
results in Simon~\cite[Chapter~4]{BS05} (see Section~\ref{sec:practests});
\item Show how two important classes of examples, systems of Schr\"odinger operators
and arbitrary order scalar operators, generate operators
$|V|^{1/2}(\pa-A_0)^{-1}U|V|^{1/2}$ which are trace class.
We reveal how the trace class properties of this Birman--Schwinger formulation
naturally reduce to the trace class properties of the 
example operators directly.
These examples also demonstrate how many practical systems 
generate such trace class operators, with the matrix trace of $V$
also equal to zero (see Section~\ref{sec:examples});
\item Prove simply and directly that the ratio of 
the Evans function and free Evans function equals the ratio of 
the transmission coefficient and free transmission coefficient. 
This new insight clarifies their relationship and indicates a 
convenient rescaling of the state variables that normalizes the 
free transmission coefficient to unity (see Section~\ref{sec:dets});
\item We show the matrix trace of the semi-separable kernels
of Birman--Schwinger operators $|V|^{1/2}(\pa-A_0)^{-1}U|V|^{1/2}$ are continuous
along the diagonal, despite the fact the kernels have a jump
discontinuity there---this assumes the matrix trace of $V$ is zero.
Hence we can unambiguously define the trace of such operators.
We then provide a simple and direct proof that the 
scaled transmission coefficient equals the Fredholm determinant of
$|V|^{1/2}(\pa-A_0)^{-1}U|V|^{1/2}$, assuming it is trace class
(see Section~\ref{sec:equivtheorem}); and
\item Prove, for the case of distinct far fields, the 
the scaled transmission coefficient equals the Fredholm determinant
with mild algebraic decay constraints on $V$
(see Section~\ref{sec:distinctfarfields}). 
\end{enumerate}
Items (ii) to (v) above are a self-contained collection of new results.
We bookend the sections above with Sections~\ref{sec:characterizations} and
\ref{sec:conclu}. In Section~\ref{sec:characterizations}
we provide preliminary results characterizing the spaces of trace class 
and Hilbert--Schmidt class operators and their relation. We include
some important inequalities required in subsequent sections. 
In Section~\ref{sec:conclu} we summarize our results, discuss
conclusions we can draw from them and outline possible future projects.

\section{\textsf{Characterizations}}\label{sec:characterizations}
To be self-contained we record a few basic facts on 
compact operators that we shall need. We refer to Reed \& Simon~\cite{RSI,RSII},
Simon~\cite{BS05} and Gohberg, Goldberg \& Krupnik~\cite{GGK}
for more details. Let $\mathbb H$ denote a 
separable Hilbert space with unitary basis $\{\varphi_m\}_{m\geqslant1}$
and standard inner product $\la\;\cdot\;,\;\cdot\;\ra_{\Hb}$.
We use $\fJ_\infty=\fJ_\infty(\Hb)$ to denote
the set of compact operators in $\Hb$. An operator
$K\in\fJ_\infty$ is positive if $\la\varphi,K\varphi\ra_{\Hb}\geqslant0$
for all $\varphi\in\Hb$. For any positive operator $K$ there is
a unique operator $\sqrt{K}$ such that $K=(\sqrt{K})^2$. The 
adjoint operator $K^\dag$ to $K$ is the unique operator such
that $\la K^\dag\varphi,\varphi\ra_{\Hb}=\la\varphi,K\varphi\ra_{\Hb}$ 
for all $\varphi\in\Hb$. The operator $K^\dag K$ is positive
as $\la K^\dag K\varphi,\varphi\ra_{\Hb}=\|K\varphi\|_{\Hb}^2\geqslant0$.
In particular we define $|K|=\sqrt{K^\dag K}$. Lastly there exists
a unique unitary operator $U$ such that $K=U\,|K|$.
For any operator $K\in\fJ_{\infty}$, we define its trace by 
$\mathrm{tr}\,K\coloneqq\sum_{m\geqslant1}\la\varphi_m,K\varphi_m\ra_{\Hb}$.
When it exists, the trace is linear and independent of the 
unitary basis chosen. 
The Schatten--von Neumann classes of compact operators $\fJ_p=\fJ_p(\Hb)$ 
for any $p\geqslant1$ are then defined as follows,
\begin{equation*}
\fJ_p\coloneqq\{K\in\fJ_\infty\colon\mathrm{tr}\,|K|^p<\infty\}.
\end{equation*}
The set $\fJ_p$ equipped with the norm $\|K\|^p_{\fJ_p}\coloneqq\mathrm{tr}\,|K|^p$ 
is a Banach space. An operator $K\in\fJ_\infty$ is \emph{trace class} if 
it belongs to $\fJ_1$ and \emph{Hilbert--Schmidt class} if it belongs to $\fJ_2$. 
The latter class $\fJ_2$ is a Hilbert space with inner product 
$\la K_1,K_2\ra_{\fJ_2}\coloneqq\mathrm{tr}K_1^\dag K_2$. 
A crucial property of the trace is that the trace of a product composition, of any 
bounded operator with a trace class operator, is invariant to their permutation. 
We can also characterize the Schatten--von Neumann classes of compact operators
$\fJ_p$ as follows. The eigenvalues $\{\lambda_m\}_{m\geqslant1}$ of any compact 
operator $K\in\fJ_\infty$ are finite in number away from the origin and the 
origin itself is the only possible accumulation point. The singular values 
$\{s_m\}_{m\geqslant1}$ of $K\in\fJ_\infty$ are the eigenvalues of $\sqrt{K^\dag K}$.
Then we can equivalently characterize 
$\mathrm{tr}\,K^p=\sum_{m\geqslant1}\lambda_m^p$ and
$\mathrm{tr}\,|K|^p=\sum_{m\geqslant1}s_m^p$.
The former is bounded by the latter.
There is a natural ordering of the Schatten--von Neumann classes as follows:
$\fJ_p\hookrightarrow\fJ_q$ for any $p\leqslant q$. 
Fundamentally, for any trace class operator $K\in\fJ_1$
the Fredholm determinant
$\mathrm{det}_1(\id+\epsilon K)\coloneqq\prod_{m\geqslant1}(1+\epsilon\lambda_m)$,
is entire in 
$\epsilon\in\C$. 
Using the relation $\mathrm{det}\,\exp=\exp\,\mathrm{tr}$
we can also characterize it (for $p=1$) by
\begin{equation*}
\mathrm{det}_p(\id+\epsilon K)
=\exp\,\sum_{\ell\geqslant p}\frac{(-1)^{\ell-1}}{\ell}
\epsilon^\ell\mathrm{tr}\,K^\ell.
\end{equation*}
When $p$ is an integer greater than one, we define the $p$-modified 
or regularized Fredholm determinants for compact operators $K\in\fJ_p$ 
by this last formula as well, knocking out the lower order non-convergent traces.
Three further results will prove very useful to us. First, if 
$A,B\in\fJ_2$ then $AB\in\fJ_1$. Second, if $A\colon\Hb\to\Hb$ 
is a bounded operator and $B\in\fJ_1$, then $AB\in\fJ_1$ and $BA\in\fJ_1$.
This is the trace class ideal property.
Third, if $A\colon\Hb\to\Hb$ is a bounded operator and $B\in\fJ_2$, 
then $AB\in\fJ_2$ and $BA\in\fJ_2$. This is the Hilbert--Schmidt ideal property.
Indeed we have, 
\begin{equation*}
\|AB\|_{\fJ_1}\leqslant\|A\|_{\fJ_2}\|B\|_{\fJ_2},\quad
\|AB\|_{\fJ_1}\leqslant\|A\|_{\mathrm{op}}\|B\|_{\fJ_1},
\quad\text{and}\quad
\|AB\|_{\fJ_2}\leqslant\|A\|_{\mathrm{op}}\|B\|_{\fJ_2},
\end{equation*}
which also hold for $BA$ and where $ \|\,\cdot\,\|_{\mathrm{op}}$ 
denotes the operator norm. The proof of these three results can be found 
for example in Conway~\cite[Section~18]{Conway}.

\section{\textsf{Practical tests}}\label{sec:practests}
The natural setting we require, and which we assume hereafter, 
is the separable Hilbert space of $\C^n$-valued
square integrable functions $\Hb=L^2(\R;\C^n)$;
see Reed \& Simon~\cite[p.~121]{RSII} for an example basis. 
Since we will be concerned with kernel functions, we
also require the separable Hilbert space $L^2(\R^2;\C^{n\times n})$
with inner product, for any $G,H\in L^2(\R^2;\C^{n\times n})$, given by
\begin{equation*}
\la G,H\ra_{L^2(\R^2;\C^{n\times n})}\coloneqq 
\int_{\R^2}\mathrm{tr}\,\bigl(G^\dag(x;y)H(x;y)\bigr)
\,\mathrm{d}x\mathrm{d}y.
\end{equation*}
The following fundamental lemma is proved in Appendix~\ref{app:HSlemma}.
%
\begin{lemma}[\textbf{\textsf{Hilbert--Schmidt class operators}}]\label{lemma:HS}
The operator $K\in\fJ_\infty$ is Hilbert--Schmidt if and only if there
is a function $G\in L^2(\R^2;\C^{n\times n})$ such that
\begin{equation*}
(K\varphi)(x)=\int_{\R}G(x;y)\,\varphi(y)\,\mathrm{d}y,
\end{equation*}
for all $\varphi\in L^2(\R;\C^n)$. In addition 
we have $\|K\|_{\fJ_2}=\|G\|_{L^2(\R^2;\C^{n\times n})}$.
\end{lemma}
The Fourier transform of operators will play a key role in our analysis.
We define the Fourier transform $L^2(\R;\C^n)\to L^2(\R;\C^n)$ 
and inverse Fourier transform as the maps
$\varphi\mapsto\hat\varphi$ and $\hat\varphi\mapsto\varphi$ respectively 
given by 
$\hat\varphi(\xi)\coloneqq(2\pi)^{-1/2}
\int_{\R}\varphi(x)\mathrm{e}^{-\mathrm{i}\xi x}\,\mathrm{d}x$ and
$\varphi(x)\coloneqq(2\pi)^{-1/2}\int_{\R}
\hat\varphi(\xi)\mathrm{e}^{\mathrm{i}x\xi}\,\mathrm{d}\xi$.
Suppose an operator $K^\ast\colon L^2(\R;\C^n)\to L^2(\R;\C^n)$ is such that
its Fourier transform $\hat K=\hat K(\xi)$ acts \emph{multiplicatively} 
in Fourier space, i.e.\/ we have 
$\widehat{K^\ast\varphi}=\hat K\hat\varphi$,
where the product $\hat K\hat\varphi$ is matrix multiplication. 
Given any multiplicative operator $\hat K=\hat K(\xi)$ in Fourier space with
integral kernel $G$ in physical space, we think of $K^\ast$ 
as the map taking $\varphi$ to the inverse Fourier transform
of $\hat K\hat\varphi$, or equivalently, 
$K^\ast\colon\varphi\mapsto(2\pi)^{-1/2}\,G\ast\varphi$,
where $G\ast\varphi$ represents the convolution of $G$ and $\varphi$.
We note that if $H\in L^2(\R;\C^{n\times n})$ then
\begin{equation*}
\bigl((K^\ast H)\varphi\bigr)(x)
=(2\pi)^{-1/2}\int_{\R}G(x-y)(H\varphi)(y)\,\mathrm{d}y,
\end{equation*}
for all $\varphi\in L^2(\R;\C^n)$. Hence the kernel of $K^\ast H$ 
is $(2\pi)^{-1/2}G(x-y)H(y)$.
We now prove the following lemma which is the matrix version 
of a result given in Simon~\cite[Chapter~4]{BS05}.
\begin{lemma}[\textbf{\textsf{Practical test for Hilbert--Schmidt class}}]
\label{lemma:HSpractest}
Suppose $\hat K,H\in L^2(\R;\C^{n\times n})$ then $K^\ast H\in\fJ_2$, 
and indeed we have 
\begin{equation*}
\|K^\ast H\|_{\fJ_2}\leqslant(2\pi)^{-1/2}
\|\hat K\|_{L^2(\R;\C^{n\times n})}\|H\|_{L^2(\R;\C^{n\times n})}.
\end{equation*}
\end{lemma}
\begin{proof}
By direct computation, line by line we successively use 
the following results: (i) The kernel of $K^\ast H$
is $(2\pi)^{-1/2}G(x-y)H(y)$ and the trace of a product of two operators 
is invariant to their permutation; (ii) 
The Cauchy--Bunyakovski--Schwarz inequality in the form 
$\mathrm{tr}\,A^\dag B\leqslant (\mathrm{tr}\,A^\dag A)^{1/2}(\mathrm{tr}\,B^\dag B)^{1/2}$
for any two matrices $A$ and $B$---see Meyer~\cite[p.~289]{Meyer}. 
We used this inequality with $A=G^\dag G$ and $B=HH^\dag$ and also that 
$\mathrm{tr}\,(HH^\dag)^\dag HH^\dag\equiv \mathrm{tr}\,(H^\dag H)^\dag H^\dag H$;
(iii) The sum of the squares of singular values
is less than the square of their sum, i.e.\/ 
$(\mathrm{tr}\,A^\dag A)^{1/2}\leqslant \mathrm{tr}\,(A^\dag A)^{1/2}$;
(iv) The Young inequality; (v) That 
$\bigl\|\mathrm{tr}\,G^\dag G\bigr\|_{L^1(\R;\C)}=\|G\|_{L^2(\R;\C^{n\times n})}^2$
and (vi) The Plancherel Theorem. The direct computation is as follows, 
\begin{align*}
\|K^\ast H\|_{\fJ_2}^2
&=\frac{1}{2\pi}\int_{\R^2}\mathrm{tr}\,\Bigl(\bigl(G^\dag G\bigr)(x-y)\,
\bigl(HH^\dag\bigr)(y)\Bigr)\,\mathrm{d}y\,\mathrm{d}x\\
&\leqslant\frac{1}{2\pi}\int_{\R^2}
\Bigl(\mathrm{tr}\,\bigl((G^\dag G)^\dag(G^\dag G)\bigr)(x-y)\Bigr)^{1/2}
\Bigl(\mathrm{tr}\,\bigl((H^\dag H)^\dag(H^\dag H)\bigr)(y)\Bigr)^{1/2}
\,\mathrm{d}y\,\mathrm{d}x\\
&\leqslant\frac{1}{2\pi}\int_{\R^2}\bigl(\mathrm{tr}\,G^\dag G\bigr)(x-y)
\cdot\bigl(\mathrm{tr}\,H^\dag H\bigr)(y)\,\mathrm{d}y\,\mathrm{d}x\\
&\leqslant\frac{1}{2\pi}\bigl\|\mathrm{tr}\,G^\dag G\bigr\|_{L^1(\R;\C)}
\bigl\|\mathrm{tr}\,H^\dag H\bigr\|_{L^1(\R;\C)}\\
&=\frac{1}{2\pi}\|G\|_{L^2(\R;\C^{n\times n})}^2\|H\|_{L^2(\R;\C^{n\times n})}^2\\
&=\frac{1}{2\pi}\|\hat K\|_{L^2(\R;\C^{n\times n})}^2\|H\|_{L^2(\R;\C^{n\times n})}^2.
\end{align*}
\qed
\end{proof}
Lemma~\ref{lemma:HSpractest} provides us with a test to determine if
a given bounded operator is of Hilbert--Schmidt class. In practice, 
suppose we know the Fourier transform $\hat K=\hat K(\xi)$ of an operator 
and we have established it lies in $L^2(\R;\C^{n\times n})$. 
Further suppose $J=J(x)$ and $H=H(x)$ are bounded multiplicative 
operators from $L^2(\R;\C^n)$ to $L^2(\R;\C^n)$ and 
$J\in L^\infty(\R;\C^{n\times n})$ and $H\in L^2(\R;\C^{n\times n})$.
The kernel of $JK^\ast H$ is $(2\pi)^{-1/2}J(x)G(x-y)H(y)$. 
The Hilbert--Schmidt ideal property 
$\|JK^\ast H\|_{\fJ_2}\leqslant\|J\|_{\mathrm{op}}\|K^\ast H\|_{\fJ_2}$
and Lemma \ref{lemma:HSpractest} reveal that 
\begin{equation*}
\|JK^\ast H\|_{\fJ_2}\leqslant(2\pi)^{-1/2}
\|J\|_{L^\infty(\R;\C^{n\times n})}
\|\hat K\|_{L^2(\R;\C^{n\times n})}\|H\|_{L^2(\R;\C^{n\times n})}.
\end{equation*}
We would like an analogous practical test of when an operator such as $JK^\ast H$
from $L^2(\R;\C^{n})$ to $L^2(\R;\C^{n})$ is of trace class. To achieve this we require the 
two classical results mentioned above, that the product of two Hilbert--Schmidt
class operators is of trace class, and the trace class ideal property.
To establish that $K^\ast H$ is trace class, for example where $H=H(x)$
is a bounded multiplicative operator and $K^\ast$ is the operator corresponding
to $\hat K=\hat K(\xi)$ in Fourier space, we naturally require the stronger
conditions $\hat K\in L_w^2(\R;\C^{n\times n})$ and $H\in L_w^2(\R;\C^{n\times n})$,
i.e.\/ in a weighted square integrable space. More precisely we define
\begin{equation*}
L_w^2(\R;\C^{n\times n})\coloneqq\bigl\{F\in L^2(\R;\C^{n\times n})
\colon Fw\in L^2(\R;\C^{n\times n})\bigr\},
\end{equation*}
where the weight function $w$ is defined via the map 
$w\colon \xi\mapsto (1+\xi^2)^{1/2}$.
The practical test takes the following form, which is the matrix version
of a result given in Reed \& Simon~\cite[Appendix~2]{RSclassicalscattering}.
%
\begin{lemma}[\textbf{\textsf{Practical test for trace class}}]\label{lemma:tracepractest}
Suppose $\hat K,H\in L_w^2(\R;\C^{n\times n})$ 
then $K^\ast H\in\fJ_1$, and there exists a constant $c>0$ such that
\begin{equation*}
\|K^\ast H\|_{\fJ_1}\leqslant 
c\,\|\hat K\|_{L_w^2(\R;\C^{n\times n})}
\|H\|_{L_w^2(\R;\C^{n\times n})}.
\end{equation*}
\end{lemma}
\begin{proof}
We adapt the proof for the scalar case given in 
Reed \& Simon~\cite[Appendix~2]{RSclassicalscattering}.
Using the practical test for Hilbert--Schmidt class
Lemma~\ref{lemma:HSpractest} and that 
$L_w^2(\R;\C^{n\times n})\hookrightarrow L^2(\R;\C^{n\times n})$, 
our assumptions on $\hat K$ and $H$ 
imply that $K^\ast H$ is of Hilbert--Schmidt class and
has integral kernel $(2\pi)^{-1/2}G(x-y)H(y)$.
We decompose $K^\ast H=AB$ into the product
of the two operators $A$ and $B$ defined as follows,
$A\coloneqq K^\ast(1-\pa^2)^{1/2}w^{-1}$ and $B\coloneqq w(1-\pa^2)^{-1/2}H$,
where $w=w(x)$ is the weight function. Then using 
$\|AB\|_{\fJ_1}\leqslant\|A\|_{\fJ_2}\|B\|_{\fJ_2}$ 
and slightly adapting the proof of Lemma~\ref{lemma:HSpractest}
to take into account that $w=w(x)$ is scalar, we find
\begin{align*}
\|K^\ast H\|_{\fJ_1}&=\bigl\|(K^\ast(1-\pa^2)^{1/2}w^{-1})(w(1-\pa^2)^{-1/2}H)\bigr\|_{\fJ_1}\\
&\leqslant\bigl\|K^\ast(1-\pa^2)^{1/2}w^{-1}\bigr\|_{\fJ_2}
\bigl\|w(1-\pa^2)^{-1/2}H\bigr\|_{\fJ_2}\\
&\leqslant(2\pi)^{-1/2}\|\hat Kw\|_{L^2(\R;\C^{n\times n})}\bigl\|w^{-1}\bigr\|_{L^2(\R;\R)}
\bigl\|w(1-\pa^2)^{-1/2}H\bigr\|_{\fJ_2}\\
&=2^{-1/2}\|\hat K\|_{L_w^2(\R;\C^{n\times n})}\bigl\|w(1-\pa^2)^{-1/2}H\bigr\|_{\fJ_2},
\end{align*}
since $\|w^{-1}\|_{L^2(\R;\R)}=\pi^{1/2}$. Note in the expression
$\hat Kw$ in the penultimate line, the term $w=w(\xi)$ represents the 
Fourier transform of $(1-\pa^2)^{1/2}$. 
The action of the operator $(1-\pa^2)^{-1/2}$ is given by
convolution with a function $G\in L^2(\R;\R)$. In addition,
since its Fourier transform $w^{-1}=w^{-1}(\xi)$ is analytic
in the strip $\{\xi\in\C\colon |\mathrm{Im}(\xi)|<1\}$, by
the Paley--Wiener Theorem (see Reed \& Simon~\cite[Theorem~IX.13]{RSII})
we also know that $G\varepsilon\in L^2(\R;\R)$, where 
$\varepsilon\colon x\mapsto\exp(|x|/2)$. The 
kernel of the operator $w(1-\pa^2)^{-1/2}H$ is 
$(2\pi)^{-1/2}w(x)\,G(x-y)\,H(y)$ and so
\begin{equation*}
\bigl\|w(1-\pa^2)^{-1/2}H\bigr\|_{\fJ_2}^2
=\frac{1}{2\pi}\int_{\R}\biggl(\int_\R (1+x^2)\,G^2(x-y)\,\mathrm{d}x\biggr)
\Bigl(\mathrm{tr}\,\bigl(H^\dag H\bigr)(y)\Bigr)\,\mathrm{d}y.
\end{equation*}
Using 
$\int_\R (1+x^2)\,G^2(x-y)\,\mathrm{d}x
\leqslant\int_\R \bigl(1+2x^2+2y^2\bigr)G^2(x)\,\mathrm{d}x\leqslant c\,(1+y^2)$
for some constant $c>0$ 
as $G, G\varepsilon\in L^2(\R;\R)$, we find 
$\bigl\|w(1-\pa^2)^{-1/2}H\bigr\|_{\fJ_2}
\leqslant c\,\|H\|_{L_w^2(\R;\C^{n\times n})}$.
We now insert this into the estimate above.
\qed
\end{proof}
The trace ideal property of $\fJ_1$ implies
$\|JK^\ast H\|_{\fJ_1}\leqslant\|J\|_{\mathrm{op}}\|K^\ast H\|_{\fJ_1}$
and thus $JK^\ast H$ is also trace class for any bounded 
multiplicative operator $J\colon L^2(\R;\C^n)\to L^2(\R;\C^n)$ if
$K^\ast H$ is. 
We end this section with a trace formula and an immediate corollary.
A proof is provided in Appendix~\ref{app:traceformula}. 
\begin{lemma}[\textbf{\textsf{Trace formula}}]\label{lemma:traceproductformula}
Given any integer $\ell\geqslant2$, suppose that $K_1$, $K_2$, \ldots, $K_\ell$ 
are Hilbert--Schmidt operators with canonical respective kernels 
$G_1$, $G_2$, \ldots, $G_\ell$ in $L^2(\R^2;\C^{n\times n})$. 
Then we have      
\begin{equation*}                                                                                 
\mathrm{tr}\,K_1K_2\cdots K_\ell=                                                                 
\mathrm{tr}\,\int_{\R^{\ell}}                                                                     
G_1(y_1;y_2)G_2(y_2;y_3)\cdots                                                                    
G_\ell(y_{\ell};y_1)\,\mathrm{d}y_{\ell}\cdots\,\mathrm{d}y_1.                                    
\end{equation*}  
If $K$ is trace class and its kernel $G\in\C^{n\times n}$
is continuous on $\R^2$, then 
$\mathrm{tr}\,K=\mathrm{tr}\,\int_{\R}G(x;x)\,\mathrm{d}x$.
\end{lemma}
\begin{corollary}[\textbf{\textsf{Trace formula: separable kernel}}]
\label{corollary:separablekernel}
If a Hilbert--Schmidt operator $K$ has a separable kernel $G(x;y)=G_1(x)G_2(y)$  
for all $(x,y)\in\R^2$, then for any integer $\ell\geqslant1$ we have  
\begin{equation*}                  
\mathrm{tr}\,K^\ell 
=\mathrm{tr}\,\Biggl(\int_\R G_2(y)G_1(y)\,\mathrm{d}y\Biggr)^\ell.                              
\end{equation*}           
\end{corollary}
\begin{remark} Some additional comments on the material above are as follows:
(i) The map $K^\ast$ above corresponds to the operator $K(-\mathrm{i}\nabla)$
in Reed \& Simon~\cite[p.~57]{RSII} and Simon~\cite[p.~37]{BS05};
(ii) The sufficiency conditions in Lemma~\ref{lemma:tracepractest}
on $\hat K$ and $H$ can be weaker, for example a result of 
Birmann--Solomajk 
requires they need only be in the space of $\ell^1$-summable 
$L_w^2(I_k;\C^{n\times n})$-norms of $\hat K$ and $H$,
where $I_k$ is the unit interval with centre 
$k\in\mathbb Z$---this space contains 
$L_w^2(\R;\C^{n\times n})$---see Simon~\cite[Chapter~4]{BS05}
for more details. However the sufficient conditions quoted 
above are adequate for the applications we have in mind; 
(iii) The proof of the trace class Lemma~\ref{lemma:tracepractest} 
for scalar operators by Reed \& Simon~\cite[Appendix~2]{RSclassicalscattering} 
is for any spatial dimension $d$, with weights $w^\alpha$ and $\alpha>d/2$;
and (iv) In the trace formula Lemma~\ref{lemma:traceproductformula}, 
we relax the continuity condition on $G$ in Section~\ref{sec:equivtheorem} 
and allow $G$ to have a jump across its diagonal---the 
formula for $\mathrm{tr}\,K$ remains meaningful as the kernels
we consider have continuous matrix traces along the 
diagonal; see also Remarks~\ref{remarks:Brislawn} where we 
discuss results by Brislawn~\cite{Br88,Br91} for discontinuous kernels.
\end{remark}

\section{\textsf{Examples}}\label{sec:examples}
Anticipating our main result in Sections~\ref{sec:equivtheorem} 
\&~\ref{sec:distinctfarfields} we will be concerned with establishing
whether operators from $L^2(\R;\C^n)$ to $L^2(\R;\C^n)$ of the form  
$|V|^{1/2}(\pa-A_0)^{-1}U|V|^{1/2}$
are of trace class. Here $A_0\in\C^{n\times n}$
is a constant matrix; $V=V(x)$ is a bounded multiplicative 
matrix operator from $L^2(\R;\C^n)$ to $L^2(\R;\C^n)$ and $U=U(x)$ is
the unitary matrix obtained from the polar decomposition of $V$.
As we saw in the last section the practical
tests for  Hilbert--Schmidt or trace class operators at our disposal rely 
on testing the integrability properties
of the multiplicative operator 
$(\mathrm{i}\xi\id-A_0)^{-1}$ in Fourier space
corresponding to the operator $(\pa-A_0)^{-1}$, as well as the integrability 
properties of $|V|^{1/2}$ and $U|V|^{1/2}$.
We assume in this section that the eigenvalues of $A_0$ are non-zero and never
purely imaginary. Hence the integrability properties of the operator
$(\mathrm{i}\xi\id-A_0)^{-1}$ rely on the rate of its asymptotic decay as $|\xi|\to\infty$. 
We observe $(\mathrm{i}\xi\id-A_0)^{-1}$ is square integrable.
We assume \emph{hereafter} that 
\begin{equation*}
V\in L_{w^2}^1(\R;\C^{n\times n})\cap L^\infty(\R;\C^{n\times n})\cap C(\R;\C^{n\times n})
\end{equation*}
where $C(\R;\C^{n\times n})$ represents the space of 
continuous $\C^{n\times n}$-valued matrix functions on $\R$. Since 
$L_{w^2}^1(\R;\C^{n\times n})\cap L^\infty(\R;\C^{n\times n})
\hookrightarrow L_w^2(\R;\C^{n\times n})$,
and $L_w^2(\R;\C^{n\times n})\hookrightarrow L^2(\R;\C^{n\times n})$
our functions $V$ are square integrable.
We observe, using the unitary properties of $U$ and that $|V|$ is selfadjoint, that 
$\bigl\|U|V|^{1/2}\bigr\|_{L^2(\R;\C^{n\times n})}^2
=\bigl\||V|\bigr\|_{L^1(\R;\C^{n\times n})}$.
We thus deduce $U|V|^{1/2}\in L^2(\R;\C^{n\times n})$ as
$L_{w^2}^1(\R;\C^{n\times n})\hookrightarrow L^1(\R;\C^{n\times n})$. 
Further $|V|^{1/2}\in L^\infty(\R;\C^{n\times n})$ by assumption.
Hence by the practical test for Hilbert--Schmidt class Lemma~\ref{lemma:HSpractest},
and the immediate discussion after, the operator $|V|^{1/2}(\pa-A_0)^{-1}U|V|^{1/2}$ 
is of Hilbert--Schmidt class.

Recall 
the sufficient conditions for 
$|V|^{1/2}(\pa-A_0)^{-1}U|V|^{1/2}$ 
to pass the practical test for trace class Lemma~\ref{lemma:tracepractest}.
These are, first that $|V|^{1/2}$ have 
bounded operator norm, 
and, second that 
$U|V|^{1/2}\in L_w^2(\R;\C^{n\times n})$ which is equivalent to 
$|V|^{1/2}\in L_w^2(\R;\C^{n\times n})$,
which is equivalent to $V\in L_{w^2}^1(\R;\C^{n\times n})$.
These two conditions are satisfied by our assumptions
on the potential function stated above. However the third
condition is not satisfied as 
$(\mathrm{i}\xi\id-A_0)^{-1}\not\in L_w^2(\R;\C^{n\times n})$.
Hence in general, without further insight and knowledge, 
the operator $|V|^{1/2}(\pa-A_0)^{-1}U|V|^{1/2}$
fails the practical test for trace class Lemma~\ref{lemma:tracepractest}.
However we now consider two examples of operators of the form above
which pass the practical test for trace class. This is achieved
by taking advantage of the structure of $(\pa-A_0)^{-1}$ and $V$,
as will be apparent.
%
%
\begin{example}[System of elliptic operators]
Consider the coupled Schr\"odinger operator of the form
$\pa^2-d_0\pa-c_0-v$,
where $c_0$ and $d_0$ are constant square matrices 
and $v$ is a matrix potential on $\R$. 
Such operators for example, arise
in the study of the linear stability of travelling wave solutions to 
systems of nonlinear reaction-diffusion equations. Linearising the 
system of equations about the travelling wave in a co-moving frame 
and assuming an exponential time dependence with growth $\lambda$ (the 
spectral parameter) generates an eigenvalue problem
of the form $(\pa^2-d_0\pa-c_0-v)\,u=\lambda a_0^{-1}u$ for the 
eigenstates $u$. Here $a_0$ is the matrix of diffusion coefficients. 
Replacing $c_0+\lambda a_0^{-1}\to c_0$, the determination of eigenvalues 
reduces to finding the zeros of the determinant of the operator above.
In phase space, the 
first order constant coefficient operator corresponding to 
$\pa^2-d_0\pa-c_0$ and matrix potential $V$ 
corresponding to $v$ have the form 
\begin{equation*}
\pa-A_0=\begin{pmatrix}\pa & -\id\\-c_0 &\pa-d_0 \end{pmatrix},
\qquad\text{and}\qquad
V=\begin{pmatrix}O & O\\ v & O \end{pmatrix},
\end{equation*}
where $O$ represents the zero matrix.
The matrices $|V|^{1/2}$ and $U$ have the form
\begin{equation*}
|V|^{1/2}=\begin{pmatrix}|v|^{1/2} & O\\ O & O \end{pmatrix}
\qquad\text{and}\qquad
\qquad U=\begin{pmatrix} O & v|v|^{-1}\\ v|v|^{-1} & O \end{pmatrix},
\end{equation*}
which can be verified by direct inspection. Additional direct
calculation reveals that
\begin{equation*}
|V|^{1/2}(\pa-A_0)^{-1}U|V|^{1/2}
=\begin{pmatrix}|v|^{1/2}[(\pa-A_0)^{-1}]_{1,2}v|v|^{-1/2} & O\\ O & O \end{pmatrix}
\end{equation*}
where $[(\pa-A_0)^{-1}]_{1,2}$ represents the top righthand square matrix block in 
$(\pa-A_0)^{-1}$ whose dimensions are half those of $(\pa-A_0)^{-1}$
itself. In Fourier space, $(\pa-A_0)^{-1}$ here has the form 
\begin{equation*}
(\mathrm{i}\xi\id-A_0)^{-1}
=\begin{pmatrix} \mathrm{i}\xi\id & -\id\\-c_0 &\mathrm{i}\xi\id-d_0 \end{pmatrix}^{-1}.
\end{equation*}
Direct computation reveals 
$[(\mathrm{i}\xi\id-A_0)^{-1}]_{1,2}
=\bigl((\mathrm{i}\xi)^2\id-(\mathrm{i}\xi)d_0-c_0\bigr)^{-1}$.
This is the operator in Fourier space corresponding to
$(\pa^2-d_0\pa-c_0)^{-1}$ as expected.
We observe it is square integrable with respect to the weight $w$ as
$\bigl((\mathrm{i}\xi)^2\id-(\mathrm{i}\xi)d_0-c_0\bigr)^{-1}(1+\xi^2)^{1/2}\sim-\xi^{-1}\id$
as $\xi\to\pm\infty$.
The trace class property of $|V|^{1/2}(\pa-A_0)^{-1}U|V|^{1/2}$
transfers to the trace class property of 
$|v|^{1/2}[(\pa-A_0)^{-1}]_{1,2}v|v|^{-1/2}$, and 
the properties we assume for $V$
transfer to $v$, and vice-versa. Hence 
by the practical test for trace class 
Lemma~\ref{lemma:tracepractest}, 
$|V|^{1/2}(\pa-A_0)^{-1}U|V|^{1/2}$ is of trace class.
\end{example}
\begin{example}[High order scalar operator]
Consider the scalar $n$th order linear operator given by
$\pa^n+a_{n-1}\pa^{n-1}+\cdots+a_1\pa+a_0+v$ where the
$a_i$, $i=0,\ldots,n-1$ are scalar constants and 
$v$ is a scalar potential function on $\R$.
If we rewrite this as
a first order system in phase space, then the 
corresponding matrix potential $V$ has the same
form as that in the last example, with the only
non-zero entry being the lower left entry which is 
$v$, however this is now the scalar entry in the
lower left position and the rest of the $n\times n$
matrix $V$ has zero entries. Similarly
the matrices $|V|^{1/2}$ and $U$ have the 
same form as in the last example, but the entries shown
are scalar with the rest of the matrix entries being zero.
Direct calculation reveals that the only non-zero entry 
in the matrix operator $|V|^{1/2}(\pa-A_0)^{-1}U|V|^{1/2}$ is the 
top left scalar entry 
$|v|^{1/2}[(\pa-A_0)^{-1}]_{1,n}v|v|^{-1/2}$,
where $[(\pa-A_0)^{-1}]_{1,n}$ represents the top righthand scalar entry  
in the operator $(\pa-A_0)^{-1}$. Direct calculation reveals 
the corresponding entry in the Fourier transform of $(\pa-A_0)^{-1}$ is given by
$[(\mathrm{i}\xi\id-A_0)^{-1}]_{1,n}=\bigl((\mathrm{i}\xi)^n\id+(\mathrm{i}\xi)^{n-1}a_{n-1}
+\cdots+(\mathrm{i}\xi)a_1+a_0\bigr)^{-1}$.
This is square integrable in Fourier space 
with respect to the weight $w$.
Hence, as in the last example, 
by the practical test for trace class 
Lemma~\ref{lemma:tracepractest}, the operator
$|V|^{1/2}(\pa-A_0)^{-1}U|V|^{1/2}$ is trace class.
\end{example}
%
%
%

\section{\textsf{Evans function and transmission coefficient}}\label{sec:dets}
Our main practical concern are eigenvalue problems associated with the
stability of travelling pulses. A wide class of such eigenvalue problems 
can be expressed in the following form on $\R$: 
\begin{equation*}
(\pa-A_0-V)Y=O.
\end{equation*}
Here $A_0=A_0(\lambda)$ is a constant $\C^{n\times n}$-valued matrix
which depends linearly on the spectral parameter $\lambda\in\C$. The 
$\C^{n\times n}$-valued matrix function $V=V(x)$ with $x\in\R$ represents
a potential perturbation. We assume throughout that $V$ is $w^2$-integrable,
uniformly bounded and continuous, exactly as outlined at
the beginning of the last Examples Section~\ref{sec:examples}. 
A comprehensive reference at this stage for 
the present material is Sandstede~\cite{BS02}.
Our goal is to determine the values of $\lambda$ for which 
the operator $\pa-A_0-V$ is not invertible, i.e.\/ the spectrum
of this operator. The complement to the spectrum in $\C$ is the 
resolvent set. Specifically from the spectrum, we are 
interested in determining the pure-point spectrum of this operator
rather than its complement the essential spectrum.
To achieve this we can construct determinant discriminants which 
are analytic in the spectral parameter away from the essential
spectrum and only zero at pure-point spectrum values. Contour
integration using such determinants in the spectral parameter plane 
then provides a global and local location strategy for the pure-point
spectrum. We will assume that away from the essential spectrum
the matrix $A_0=A_0(\lambda)$ is strictly hyperbolic; this is characteristic
of a wide class of travelling pulse stability problems. 
Hence away from the essential spectrum the eigenvalue equation
above has an exponential dichotomy and the unbounded operator $\pa-A_0-V$ 
is Fredholm with index zero. At pure-point spectrum values the kernel
of $\pa-A_0-V$ is non-trivial. The existence of the exponential
dichotomy to the eigenvalue equation above, and thus Fredholm
property of $\pa-A_0-V$, as well as the locale 
of the essential spectrum, is determined by the classification of 
the solutions to the following associated constant coefficient equation which
does not exhibit a pure-point spectrum:
\begin{equation*}
(\pa-A_0)Y_0=O.
\end{equation*}
The existence of an exponential dichotomy to the eigenvalue equation 
(with potential $V$) implies that there is a, say, $k$-dimensional subspace
of solutions to the eigenvalue equation that decays exponentially
as $x\to-\infty$, and an $(n-k)$-dimensional subspace of solutions
that decays exponentially as $x\to+\infty$. We denote by $Y^-=Y^-(x;\lambda)$
the $\C^{n\times k}$-valued function whose column span coincides with
the subspace decaying as $x\to-\infty$, and by $Y^+=Y^+(x;\lambda)$
the $\C^{(n-k)\times k}$-valued function whose column span coincides with
the subspace decaying as $x\to+\infty$. Corresponding subspaces
of commensurate dimension exist for the constant coefficient equation,
we denote by $Y_0^\pm=Y_0^\pm(x;\lambda)$ the corresponding matrix
valued functions whose columns span the respective exponentially 
decaying subspaces. Asymptotically as $x\to\pm\infty$ we have
$Y^\pm\sim Y_0^\pm$.

From this perspective, the pure-point spectrum is determined by the
values of $\lambda\in\C$ for which the two subspaces spanned by the columns
of $Y^-$ and $Y^+$, intersect. In other words, a solution to the eigenvalue
problem exists that decays exponentially to zero in both far-fields.
This is equivalent to the condition that the columns of $Y^-$ and $Y^+$
are not linearly independent, and a test for that is whether the
determinant, of the $n\times n$ matrix whose columns are the columns
of $Y^-$ and $Y^+$, is zero. This intersection property should be
$x$-independent and by Liouville's Theorem an appropriate scalar factor 
achieves this. This is the Evans function, 
first introduced by Evans~\cite{JE75}. A comprehensive 
study is provided in Alexander, Gardner \& Jones~\cite{AGJ90}.
\begin{definition}[\textbf{\textsf{Evans function}}]
The Evans function is the $\lambda$-dependent complex scalar quantity
\begin{equation*} 
\exp\biggl(-\mathrm{tr}\,A_0(\lambda)x-\int_0^x\mathrm{tr}\,V(y)\,\mathrm{d}y\biggr)\,
\mathrm{det}\begin{pmatrix} Y^-(x;\lambda) & Y^+(x;\lambda)\end{pmatrix}.
\end{equation*}
%
The free Evans function is that corresponding to the 
constant coefficient operator $\pa-A_0$, i.e.\/ corresponding to the free state. 
If we replace $Y^\pm$ by $Y_0^\pm$ and set $V\equiv 0$ in the Evans function,
we get 
$\exp\bigl(-\mathrm{tr}\,A_0(\lambda)x\bigr)\,
\mathrm{det}\begin{pmatrix}Y_0^-(x;\lambda) & Y_0^+(x;\lambda)\end{pmatrix}$.
\end{definition}
Associated with the eigenvalue problem above is the adjoint 
eigenvalue problem given by $(\pa+A_0^\dag+V^\dag)Z^\dag=O$ or
$\pa Z+ZA_0+ZV=O$,
where $Z^\dag$ is $\C^{n\times 1}$-valued and $Z$ is $\C^{1\times n}$-valued.
The adjoint operator in $L^2(\R;\C^n)$ to $\pa-A_0-V$ 
is $-\pa-A_0^\dag-V^\dag$. We denote by $Z^-=Z^-(x;\lambda)$
and $Z^+=Z^+(x;\lambda)$ the, respectively, $\C^{(n-k)\times n}$-valued and
$\C^{k\times n}$-valued functions whose respective row spans determine
solution subspaces to the adjoint eigenvalue problem that
decay exponentially as $x\to-\infty$ and $x\to+\infty$. 
In addition we denote by $Z_0^\pm=Z_0^\pm(x;\lambda)$ the 
corresponding matrices for the adjoint constant coefficient equation
\begin{equation*}
\pa Z_0+Z_0A_0=O.
\end{equation*}
The solutions $Y_0^\pm$ and $Z_0^\pm$ to the constant coefficient problems 
above satisfy a diagonal relation that will be helpful in out subsequent analysis.
\begin{lemma}[\textbf{\textsf{Diagonal relation}}]\label{lemma:diagrelation}
The solutions $Y_0^\pm$ and $Z_0^\pm$ satisfy the relation
\begin{equation*}
\begin{pmatrix}
Z_0^+\\ Z_0^-
\end{pmatrix}
\begin{pmatrix}
Y_0^- & Y_0^+
\end{pmatrix}=D,
\end{equation*}
where $D$ is a constant diagonal matrix with non-zero entries.
\end{lemma}
\begin{proof}
Consider any pair of solutions $Y_0\in\C^{n\times1}$ and $Z_0\in\C^{1\times n}$ 
to the their respective constant coefficient problems. Then we see
$\pa(Z_0\,Y_0)=-Z_0A_0Y_0+Z_0A_0Y_0=O$. Thus $Z_0\,Y_0$ is constant on $\R$. 
Each solution $Y_0$ has the form $U\exp(\mu x)$ where $\mu$ is an
eigenvalue and $U\in\C^{n\times1}$ a corresponding right eigenvector of $A_0$. By our 
strict hyperbolicity assumption, there are $n$ independent solutions,
and $k$ of the eigenvalues have a positive real part.  
Each solution $Z_0$, corresponding to an adjoint solution, has the  
form $W\exp(-\nu x)$ where $\nu$ is an eigenvalue and $W\in\C^{1\times n}$ 
a corresponding left eigenvector of $A_0$. Classically if $\mu\neq\nu$
then $WU=0$, while if $\mu=\nu$ we have $WU\neq0$; 
see Meyer~\cite[p.~405, 523]{Meyer}.
\qed
\end{proof}
With all this in hand, we can now motivate and define the transmission
coefficient. Starting with the Evans function, we find 
\begin{align*}
\exp&\biggl(-\mathrm{tr}\,A_0(\lambda)x-\int_0^x\mathrm{tr}\,V(y)\,\mathrm{d}y\biggr)
\,\mathrm{det}\begin{pmatrix}Y^-& Y^+\end{pmatrix}\\
&=\exp\biggl(-\mathrm{tr}\,A_0(\lambda)x
-\int_0^x\mathrm{tr}\,V(y)\,\mathrm{d}y\biggr)\,\frac{\mathrm{det}\Biggl(
\begin{pmatrix}
Y_0^- & Y_0^+
\end{pmatrix}
\begin{pmatrix}
Z_0^+\\ Z_0^-
\end{pmatrix}\Biggr)}
{\mathrm{det}\Biggl(\begin{pmatrix}
Z_0^+\\ Z_0^-
\end{pmatrix}
\begin{pmatrix}
Y_0^- & Y_0^+
\end{pmatrix}\Biggr)}\,
\mathrm{det}\begin{pmatrix}Y^-& Y^+\end{pmatrix}\\
&=\exp\bigl(-\mathrm{tr}\,A_0(\lambda)x\bigr)\,
\mathrm{det}\begin{pmatrix}
Y_0^- & Y_0^+
\end{pmatrix}
\exp\biggl(-\int_0^x\mathrm{tr}\,V(y)\,\mathrm{d}y\biggr)\,
\frac{\mathrm{det}
\begin{pmatrix}
Z_0^+Y^- & Z_0^+Y^+\\
Z_0^-Y^- & Z_0^-Y^+
\end{pmatrix}}
{\mathrm{det}
\begin{pmatrix}
Z_0^+Y_0^- & Z_0^+Y_0^+\\
Z_0^-Y_0^- & Z_0^-Y_0^+
\end{pmatrix}}.
\end{align*}
Note that the first two factors on the right constitute
the free Evans function which is $x$-independent and $\lambda$-dependent.
We could define a generalized transmission coefficient as the product 
of the middle exponential term and the numerator in the ratio. The 
denominator would correspond to a free generalized transmission coefficient.
However, classically, we take the limit as $x\to+\infty$ in these
latter terms. The diagonal relation in Lemma~\ref{lemma:diagrelation} 
in particular implies $Z_0^-Y_0^-=O$ and $Z_0^+Y_0^+=O$, and 
thus also that $\lim_{x\to+\infty}Z_0^+Y^+=O$.
Hence the determinants in the ratio above collapse as $x\to+\infty$ 
(using that the limit and determinant operations commute).
Cancelling off the determinant of $Z_0^-Y_0^+$ which appears 
in both the numerator and denominator we find 
\begin{multline*}
\exp\biggl(-\mathrm{tr}\,A_0(\lambda)x-\int_0^x\mathrm{tr}\,V(y)\,\mathrm{d}y\biggr)\,
\mathrm{det}\begin{pmatrix}Y^-& Y^+\end{pmatrix}\\
\equiv\exp\bigl(-\mathrm{tr}\,A_0(\lambda)x\bigr)
\mathrm{det}\begin{pmatrix}
Y_0^- & Y_0^+
\end{pmatrix}\cdot
\lim_{x\to+\infty}\exp\biggl(-\int_0^x\mathrm{tr}\,V(y)\,\mathrm{d}y\biggr)\,
\frac{\mathrm{det}\bigl((Z_0^+Y^-)(x;\lambda)\bigr)}
{\mathrm{det}\bigl((Z_0^+Y_0^-)(\lambda)\bigr)}.
\end{multline*}
Note by the diagonal relation in Lemma~\ref{lemma:diagrelation}
the quantity $Z_0^+\, Y_0^-$ is constant.
For completeness we now define the transmission coefficient and
free transmission coefficient.
\begin{definition}[\textbf{\textsf{Transmission coefficient}}]
This is defined as the $\lambda$-dependent complex scalar quantity
\begin{equation*} 
\lim_{x\to+\infty}\exp\biggl(-\int_0^x\mathrm{tr}\,V(y)\,\mathrm{d}y\biggr)
\mathrm{det}\bigl((Z_0^+Y^-)(x;\lambda)\bigr).
\end{equation*}
The free transmission coefficient is simply the quantity 
$\mathrm{det}\bigl(Z_0^+Y_0^-(\lambda)\bigr)$.
\end{definition}
The relation we derived above establishes that away from the 
essential spectrum the Evans function can be decomposed as
a product of the free Evans function and a ratio of the 
transmission coefficient and free transmission coefficient.
In other words, schematically, we have  
\begin{equation*}
\frac{\text{Evans function}}{\text{Free Evans function}}
=\frac{\text{Transmission coefficient}}{\text{Free transmission coefficient}}.
\end{equation*}
We end this section with three important observations.
First, the solutions $Y_0^-\in\C^{n\times k}$, $Y_0^+\in\C^{n\times(n-k)}$,
$Z_0^-\in\C^{k\times n}$ and $Z_0^+\in\C^{(n-k)\times n}$ to the constant coefficient 
problems above have the following explicit special forms. Let $U^-\in\C^{n\times k}$ 
and $U^+\in\C^{n\times(n-k)}$ denote the matrices whose columns are
the eigenvectors of $A_0$ respectively corresponding to the 
$k$ eigenvalues with positive real part and the $n-k$ eigenvalues
with negative real part. Also let $W^-\in\C^{n\times k}$ 
and $W^+\in\C^{n\times(n-k)}$ denote the matrices whose rows are
the left eigenvectors of $A_0$ respectively corresponding to the 
$k$ eigenvalues with negative real part and the $n-k$ eigenvalues 
with positive real part. Then we have the identifications
\begin{align*}
\begin{pmatrix}
Y_0^- & Y_0^+
\end{pmatrix}&\equiv\exp(A_0x)\begin{pmatrix}
U^- & U^+
\end{pmatrix}
\equiv\begin{pmatrix}
U^-\exp(\Lambda^-x) & U^+\exp(\Lambda^+x)
\end{pmatrix}
\intertext{and}
\begin{pmatrix}
Z_0^+\\ Z_0^-
\end{pmatrix}
&\equiv\begin{pmatrix}
W^+\\ W^-
\end{pmatrix}
\exp(-A_0x)
\equiv\begin{pmatrix}
\exp(-\Lambda^-x)W^+\\ \exp(-\Lambda^+x)W^-
\end{pmatrix},
\end{align*}
where $\Lambda^-$ and $\Lambda^+$ denote the diagonal matrices 
of the eigenvalues of $A_0$ with positive and negative real parts, 
respectively. 

Second, we can rescale the free state solutions 
$Y_0^\pm$ and $Z_0^\pm$ to the constant coefficient problems so 
they satisfy some unitary relations that are helpful for our 
subsequent analysis. We choose to rescale the adjoint solutions
$Z_0^\pm$ by rescaling $W^\pm$ as follows. 
\begin{definition}[\textbf{\textsf{Unitarily scaled solutions}}]
Suppose $D$ is the constant diagonal matrix from the 
diagonal relations Lemma~\ref{lemma:diagrelation}.
Let $D_-$ denote the upper $k\times k$ block of $D$ and 
$D_+$ the lower $(n-k)\times(n-k)$ block. We define
rescaled solutions $\hat W^\pm$ and correspondingly
$\hat Z_0^\pm$ by
\begin{equation*}
\hat W^\pm\coloneqq D_\mp^{-1}W^\pm
\qquad\text{and}\qquad
\hat Z_0^\pm\coloneqq D_\mp^{-1}Z_0^\pm.
\end{equation*}
\end{definition}
Note the solutions $\hat Z_0^\pm$ have the same exponential
form as $Z_0^\pm$, but with $W^\pm$ replaced by $\hat W^\pm$.
The nomination of $\hat Z_0^\pm$ as unitarily scaled solutions is justified
in the following.
\begin{lemma}[\textbf{\textsf{Unitary relations}}]\label{lemma:orthonormality}
The solutions $Y_0^\pm$ and $\hat Z_0^\pm$ satisfy the unitary relations
\begin{equation*}
\begin{pmatrix}
\hat Z_0^+\\ \hat Z_0^-
\end{pmatrix}
\begin{pmatrix}
Y_0^- & Y_0^+
\end{pmatrix}=\mathrm{id}
\qquad\text{and}\qquad
\begin{pmatrix}
Y_0^- & Y_0^+
\end{pmatrix}
\begin{pmatrix}
\hat Z_0^+\\ \hat Z_0^-
\end{pmatrix}=\mathrm{id}.
\end{equation*}
These are generated by corresponding unitary relations satisfied by
$U^\pm$ and $\hat W^\pm$.
\end{lemma}
\begin{proof}
From the proof of the diagonal relation Lemma~\ref{lemma:diagrelation},
we already know that 
\begin{equation*}
\begin{pmatrix}
W^+\\ W^-
\end{pmatrix}
\begin{pmatrix}
U^- & U^+
\end{pmatrix}=\begin{pmatrix} D_- & O\\O & D_+\end{pmatrix}.
\end{equation*}
Substituting $W^\pm=D_\mp\hat W^\pm$ into this relation generates
the corresponding unitary relation for $\hat W^\pm$ and $U^\pm$,
i.e.\/ with the identity matrix on the right. Hence
the inverse of the $n\times n$ matrix $\bigl(\,U^-\,\,U^+\,\bigr)$ 
exists and is given by the corresponding $n\times n$ matrix with 
$\hat W^+$ and $\hat W^-$ in the upper and lower block positions, respectively, 
as shown above. The second unitary relation, with the order of the two
block matrices with $\hat W^\pm$ and $U^\pm$ swapped round, now follows;
see Meyer~\cite[p.~117]{Meyer}.
The unitary relations for $Y_0^\pm$ and $\hat Z_0^\pm$ now
follow using their explicit exponential forms. 
\qed
\end{proof}
Third, in the next section, the unitarily scaled solutions 
$\hat Z_0^\pm$ are the natural choice for constructing the Green's
kernel associated with $(\pa-A_0)^{-1}$. They are also natural
for establishing the equivalence between the transmission coefficient and 
corresponding Fredholm determinant, as the corresponding free transmission
coefficient is unity with these scaled solutions. Any result we prove
using the unitarily scaled solutions we can recover for the original
solutions $Z_0^\pm$ by substituting the relation between the two. 
This is important, as an edifying feature of the Evans function is
that it is analytic in $\lambda$ away from the essential spectrum
and its zeros correspond to eigenvalues of $\pa-A_0-V$ with
coincident multitude. The solutions $Y_0^\pm$ and $Z_0^\pm$ can
be chosen to be analytic in $\lambda$ from the start. 
The Evans function and free Evans function are 
independent of the unitary rescaling as they are defined only
using $Y^\pm$ and $Y_0^\pm$, respectively. And, as we
should expect, the ratio of the transmission coefficient 
and free transmission coefficient is also invariant to
the unitary rescaling. We return to
these points in our conclusions Section~\ref{sec:conclu}.
\begin{remark}
Classically the transmission coefficient is defined as follows. 
Consider the solu-tions $Y^-\sim Y^-_0$ as $x\to-\infty$ above. 
These are Jost solutions. 
Away from the essential spectrum, asymptotically as $x\to+\infty$ 
we have $Y^-\sim Y_0^-a+Y_0^+b$ where the constant $k\times k$ and 
$(n-k)\times(n-k)$ matrices $a=a(\lambda)$ and $b=b(\lambda)$ are 
the transmission and reflection matrix coefficients, respectively. 
For an element of the span of $Y^-$ to be an eigenfunction, 
we require it to be asymptotically in the span of $Y_0^+b$ 
or equivalently in the span of $Y_0^+$, as $x\to+\infty$. 
Or equivalently in this limit, we require an element of the 
span of $Y^-$ to be orthogonal to the subspace of $\mathbb C^n$
that is orthogonal to the subspace spanned by the columns of $Y_0^+$.
In other words we require an element of the span of $Y^-$
to be orthogonal to the subspace spanned by the rows of $Z_0^+$.
The existence of a non-trivial linear combination of columns
of $Y^-$ to be orthogonal to each row of $Z_0^+$ amounts to 
requiring $\det(Z_0^+Y^-)$ to be zero in the limit. Modulo
the constant non-zero exponential factor we define the 
transmission coefficient as this determinant.
However we note, using the diagonal relation Lemma~\ref{lemma:diagrelation}
we find $Z_0^+Y^-\sim (Z_0^+Y_0^-) a$. In other words 
$\det\bigl(a(\lambda)\bigr)$ 
equals the ratio of the transmission coefficient to
the free coefficient above. For the classical example of the 
transmission coefficient for the scalar Schr\"odinger
operator see Kapitula \& Sandstede~\cite{KS04}, for which 
$\mathrm{tr}\,V\equiv 0$. They show in this case, 
the connection between the Evans function, transmission coefficient and
Fredholm determinant (coming next). Bridges \& Derks~\cite{BD99}
show the transmission coefficient equals the Evans function up
to a non-zero analytic factor.
\end{remark}

\section{\textsf{Equivalence theorem}}\label{sec:equivtheorem}
Our goal in this section is to show that the transmission
coefficient equals the Fredholm determinant for trace class
operators, associated with eigenvalue problems on $\R$
of the form  
\begin{equation*}
(\pa-A_0-V)Y=O.
\end{equation*}
The setting is precisely that outlined in the last section.
Let us now specify which Fredholm determinant we mean. 
We have already seen that away from the essential spectrum
$(\pa-A_0)^{-1}$ exists. 
Hence our eigenvalue problem can
be expressed in the form 
\begin{equation*}
\bigl(\id-(\pa-A_0)^{-1}V\bigr)Y=O.
\end{equation*}
As in Sections~\ref{sec:examples} \& \ref{sec:dets}, 
we assume throughout the potential perturbation $V$ is $w^2$-integrable,
uniformly bounded and continuous on $\R$.
An equivalent formulation is that of Birman--Schwinger.
If use the polar decomposition for $V=U|V|$ and set $\phi\coloneqq|V|^{1/2}Y$,
then we see 
\begin{align*}
&&\bigl(\id-(\pa-A_0)^{-1}V\bigr)Y&=O\\
\Leftrightarrow&&
\bigl(\id-(\pa-A_0)^{-1}U|V|^{1/2}|V|^{1/2}\bigr)Y&=O\\
\Leftrightarrow&&
\bigl(|V|^{1/2}-|V|^{1/2}(\pa-A_0)^{-1}U|V|^{1/2}|V|^{1/2}\bigr)Y&=O\\
\Leftrightarrow&&
\bigl(\id-|V|^{1/2}(\pa-A_0)^{-1}U|V|^{1/2}\bigr)\phi&=O.
\end{align*}
Here $|V|^{1/2}(\pa-A_0)^{-1}U|V|^{1/2}$ is the Birman--Schwinger
operator as considered in the examples Section~\ref{sec:examples}.
To proceed we establish the Green's kernel
corresponding to $|V|^{1/2}(\pa-A_0)^{-1}U|V|^{1/2}$. 
Suppose $(\pa-A_0)^{-1}$ has a representation
\begin{equation*}
(\pa-A_0)^{-1}\colon \varphi(x)\mapsto \int_\R G(x;y)\varphi(y)\,\mathrm{d}y,
\end{equation*}
for some function $G\in L^2(\R^2;\C^{n\times n})$, which in this context here
with out loss of generality, we will also assume is continuously differentiable 
everywhere except along the diagonal $y=x$. Note that $G$ also depends on $\lambda$
through $A_0=A_0(\lambda)$. However we suppress the dependence on $\lambda$
in all the relevant variables for the moment.
Classical theory implies that we require
$G=G(x;y)$ to satisfy the pair of differential equations
$\pa_xG-A_0G=\delta(x-y)\id$ and $-\pa_yG-GA_0=\delta(x-y)\id$.
Two formal calculations that can retrospectively be made rigorous
reveal this is the correct prescription
for $G=G(x;y)$. First, we observe 
$(\pa_x-A_0)\int G(x;y)\varphi(y)\,\mathrm{d}y
=\int\delta(x-y)\varphi(y)\,\mathrm{d}y=\varphi(x)$.
Second, we observe 
$\int G(x;y)(\pa_y-A_0)\varphi(y)\,\mathrm{d}y
=\int\bigl(-\pa_yG(x;y)-G(x;y)A_0)\bigr)\varphi(y)\,\mathrm{d}y
=\varphi(x)$.
We know the solutions to the equations for the Green's function
away from the diagonal already. Indeed the Green's function 
$G=G(x;y)$ with the correct decay properties in the 
far field is given by the following semi-separable form
which can be confirmed by direct substitution,
\begin{equation*}
G(x;y)\coloneqq\begin{cases} 
-Y^-_0(x)\hat Z^+_0(y),&\qquad x\leqslant y,\\
Y^+_0(x)\hat Z^-_0(y),&\qquad x>y.
\end{cases}
\end{equation*}
The second relation of the unitary relations 
in Lemma~\ref{lemma:orthonormality} is 
equivalent to $Y_0^-\hat Z_0^++Y_0^+\hat Z_0^-=\id$. 
The unitarily scaled solutions $\hat Z_0^\pm$ in the
Green's kernel thus guarantee the natural jump condition 
$G(x^+;x)-G(x^-;x)=\id$ across the diagonal $y=x$ is satisfied.
\begin{remark}\label{remark:J}
The solutions $Y_0^\pm$ and $\hat Z_0^\pm$ have simple 
exponential forms as indicated at the end of the last section. 
Indeed utilizing these, we can express $G$ in the form $G=G(x-y)$, consistent
with the kernel result for Hilbert--Schmidt operators stated 
prior to the practical test for Hilbert--Schmidt class Lemma~\ref{lemma:HSpractest}.
Indeed we note that 
$Y^\pm_0(x)\hat Z^\mp_0(y)\equiv U^\pm\exp\bigl(\Lambda^\pm(x-y)\bigr)\hat W^\mp$.
\end{remark}
The Green's kernel corresponding to $|V|^{1/2}(\pa-A_0)^{-1}U|V|^{1/2}$ 
has the semi-separable form $|V(x)|^{1/2}G(x;y)U(y)|V(y)|^{1/2}$.
Upon closer inspection the unitary relation given by
$Y_0^-\hat Z_0^++Y_0^+\hat Z_0^-=\id$ 
implies that the diagonal elements of the Green's kernel matrix 
$G=G(x;y)$ have a unit jump at the diagonal $y=x$, while the off-diagonal
elements are continuous. In addition the unitary relation implies
that along the diagonal $y=x$ we have
\begin{align*}
&&|V(x)|^{1/2}Y_0^-\hat Z_0^+U(x)|V(x)|^{1/2}
+|V(x)|^{1/2}Y_0^+\hat Z_0^-U(x)|V(x)|^{1/2}&=V(x)\\
\Rightarrow&&\mathrm{tr}\,|V(x)|^{1/2}Y_0^-\hat Z_0^+U(x)|V(x)|^{1/2}
+\mathrm{tr}\,|V(x)|^{1/2}Y_0^+\hat Z_0^-U(x)|V(x)|^{1/2}&=\mathrm{tr}\,V(x).
\end{align*}
Hence if assume $\mathrm{tr}\,V\equiv0$, then along the diagonal $y=x$ we have
\begin{equation*}
-\mathrm{tr}\,|V(x)|^{1/2}Y_0^-\hat Z_0^+U(x)|V(x)|^{1/2}
=\mathrm{tr}\,|V(x)|^{1/2}Y_0^+\hat Z_0^-U(x)|V(x)|^{1/2},
\end{equation*}
and thus the matrix trace of the kernel $|V(x)|^{1/2}G(x;y)U(y)|V(y)|^{1/2}$
is continuous at the diagonal $y=x$. We can thus unambiguously define the trace
of $|V|^{1/2}(\pa-A_0)^{-1}U|V|^{1/2}$. 
\begin{definition}[\textbf{\textsf{Trace for discontinuous kernels}}]
For potential peturbations $V$ which are $w^2$-integra-ble, uniformly
bounded and continuous on $\R$ and with $\mathrm{tr}\,V\equiv0$,
we define the trace of the linear operator 
$|V|^{1/2}(\pa-A_0)^{-1}U|V|^{1/2}$, for which $(\pa-A_0)^{-1}$
has a kernel with a jump along the diagonal, by
\begin{equation*}
\mathrm{tr}\,\bigl(|V|^{1/2}(\pa-A_0)^{-1}U|V|^{1/2}\bigr)\coloneqq
\mathrm{tr}\,\int_\R \bigl(-|V(x)|^{1/2}
Y_0^-(x)\hat Z_0^+(x)U(x)|V(x)|^{1/2}\bigr)\,\mathrm{d}x.
\end{equation*}
\end{definition}
\begin{remark}\label{remarks:Brislawn}
We note the following: (i) Consider the examples of 
trace class operators in Section~\ref{sec:examples}. 
The Green's kernels 
corresponding to 
$|V|^{1/2}(\pa-A_0)^{-1}U|V|^{1/2}$ are given by
$|v(x)|^{1/2}G_{12}(x;y)v(y)|v(y)|^{-1/2}$ and 
$|v(x)|^{1/2}G_{1n}(x;y)v(y)|v(y)|^{-1/2}$ in the first
and second examples, respectively. They both involve 
the off-diagonal elements of $G$ only and hence these
kernels are continuous; (ii) By analogous arguments to those above, 
if $\mathrm{tr}\,V\equiv0$, then along the diagonal $y=x$ the
matrix trace of the kernel $G(x;y)V(y)$ corresponding to
$(\pa-A_0)^{-1}V$ is also continuous.
Indeed, using the invariance of the trace to the order
of the product of two operators, the matrix traces of $G(x;y)V(y)$ and 
$|V(x)|^{1/2}G(x;y)U(y)|V(y)|^{1/2}$ are equal along the diagonal $y=x$;
(iii) Consider the scalar Schr\"odinger operator $\pa^2-c_0-v$ 
which is a special case of both examples in Section~\ref{sec:examples}.
Re-writing it as a first order operator and focusing on 
$(\pa-A_0)^{-1}V$, we observe if $G=G(x-y)$
is the $2\times2$ Green's kernel corresponding to $(\pa-A_0)^{-1}$,
then the Green's kernel corresponding to $(\pa-A_0)^{-1}V$
is the $2\times2$ matrix whose left column contains
$G_{12}v$ and $G_{22}v$ and whose right column is 
zero. The Fourier transforms of $G_{12}$ and $G_{22}$ are 
$\hat G_{12}(\xi)=-(\xi^2-c_0)^{-1}$ and 
$\hat G_{22}=-\mathrm{i}\xi(\xi^2-c_0)^{-1}$.
We note $\hat G_{12}$ is square integrable with respect to the
weight function $w$, but $\hat G_{12}$ is not. Indeed they are explicitly 
given by $G_{12}(x)=(\pi/2c_0)^{1/2}\exp(-c_0^{1/2}|x|)$
and $G_{22}(x)=(\pi/2c_0)^{1/2}\pa_x\exp(-c_0^{1/2}|x|)$. 
From Gohberg \textit{et al.\/ }~\cite[p.~244]{GGK} we have the following
result for scalar separable integral operators whose kernel functions either
side of the diagonal $y=x$ are continuous, up to and including the 
diagonal. If the operator is trace class then the kernel function is
continuous on $\R^2$. Our results in (i) above for the corresponding 
Birman--Schwinger operator are consistent with this. Further, since the 
scalar operator corresponding to $G_{22}$ has a jump along the diagonal $y=x$,
we conclude it cannot be trace class; and
(iv) Brislawn~\cite{Br88,Br91} has shown how to define the trace of 
a trace class operator with a kernel that is only square integrable
and thus not necessarily continuous. This is achieved
by averaging on cubes via the Hardy--Littlewood maximal function.
In particular the plain Volterra operator on $L^2([0,1];\R)$
with kernel equal to $1$ below the diagonal and $0$ above, 
has trace equal to $1/2$. However it has singular values given by
$2/(\pi(2n+1))$ and thus is not trace class; 
see Brislawn~\cite[Example~3.2]{Br88}.
\end{remark}
%
%
We now establish the main result of this section, our equivalence theorem. 
Before we state and prove it, we require the following key lemma. 
It concerns the solutions $Y^-$ of the eigenvalue
problem $(\pa-A_0-V)Y=O$ whose column span coincides with
the subspace of solutions decaying as $x\to-\infty$. Indeed 
the crucial insight is that by the variation of constants formula 
the solutions $Y^-$ satisfy the Volterra integral equation 
\begin{equation*}
Y^-=Y^-_0+\int_{-\infty}^x \exp\bigl(A_0(x-y)\bigr)V(y)Y^-(y)\,\mathrm{d}y.
\end{equation*}
\begin{lemma}[\textbf{\textsf{Volterra integral equation}}]\label{lemma:VIE}
Assume $\mathrm{tr}\,V\equiv 0$ and $A_0$ and $V=V(x)$ 
satisfy the conditions stated above.
Let $J$ denote the integral operator with kernel 
$H(x-y)|V(x)|^{1/2}\exp\bigl(A_0(x-y)\bigr)U(y)|V(y)|^{1/2}$, 
where $H=H(x)$ is the Heaviside step function. First, 
if we set $\phi^-\coloneqq|V|^{1/2}Y^-$ and 
$\phi^-_0\coloneqq|V|^{1/2}Y^-_0$, then $Y^-$ solves the 
eigenvalue problem $(\pa-A_0-V)Y=O$ with $Y^-\sim Y^-_0$ 
as $x\to-\infty$ if and only if $\phi^-$ satisfies
\begin{equation*}
\phi^-=\phi_0^-+J\phi^-.
\end{equation*}
Equivalently we have $\phi_0^-=(\mathrm{id}-J)\phi^-$.
Second, an equivalent expression for the kernel for $J$ is given by
$H(x-y)\bigl(\phi_0^-(x)\hat Z_0^+(y)+\phi_0^+(x)\hat Z_0^-(y)\bigr)U(y)|V(y)|^{1/2}$,
where $\phi^+_0\coloneqq|V|^{1/2}Y^+_0$.
\end{lemma}
\begin{proof}
Preceding the lemma, we already established that 
$Y^-$ satisfies the Volterra integral equation
shown if and only if $Y^-$ solves the 
eigenvalue problem $(\pa-A_0-V)Y=O$ with $Y^-\sim Y^-_0$ 
as $x\to-\infty$. To show $Y^-$ satisfies the Volterra 
integral equation preceding the lemma if and only if 
$\phi^-$ satisfies the Volterra integral equation shown in
the lemma, follows by premultiplying the Volterra integral
equation for $Y^-$ by $|V(x)|^{1/2}$ and using the definition
of $\phi^-$. The equivalence of the kernels follows 
from the unitary relation $U^-\hat W^++U^+\hat W^-=\id$ and that
$Y^\pm_0(x)=\exp(A_0x)U^\pm$ and $\hat Z^\pm_0(y)=\hat W^\pm\exp(-A_0y)$.
\qed
\end{proof}
\begin{theorem}[\textbf{\textsf{Equivalence}}]\label{theorem:equivalence}
Assume the operator $|V|^{1/2}(\pa-A_0)^{-1}U|V|^{1/2}$ is trace class,
$\mathrm{tr}\,V\equiv 0$ and $A_0$ and $V=V(x)$ satisfy the conditions
stated above. Then the transmission coefficient for 
the eigenvalue problem $(\pa-A_0-V)Y=O$ and Fredholm determinant 
of $\mathrm{id}-|V|^{1/2}(\pa-A_0)^{-1}U|V|^{1/2}$ are equal:
\begin{equation*}
\lim_{x\to+\infty}\mathrm{det}\bigl((\hat Z_0^+Y^-)(x)\bigr)
=\mathrm{det}_1\bigl(\mathrm{id}-|V|^{1/2}(\pa-A_0)^{-1}U|V|^{1/2}\bigr).
\end{equation*}
\end{theorem}
\begin{proof}
First we focus on the transmission coefficient.
Since we have $\pa Y^-=(A_0+V)Y^-$ and $\pa\hat Z_0^+=-\hat Z_0^+A_0$,
the product rule implies $\pa(\hat Z_0^+Y^-)=\hat Z_0^+VY^-$.
Using the unitary relations in Lemma~\ref{lemma:orthonormality} 
we know $\lim_{x\to-\infty}(\hat Z_0^+Y^-)(x)=\id$ and hence 
\begin{equation*}
\lim_{x\to+\infty}\mathrm{det}\bigl((\hat Z_0^+Y^-)(x)\bigr)
=\mathrm{det}\biggl(\id+\int_\R(\hat Z_0^+VY^-)(y)\,\mathrm{d}y\biggr).
\end{equation*}
We can also establish this result by premultiplying 
the integral equation for $Y^-$ by $\hat Z_0^+$, using the 
unitary relations, and then applying the commutable operations
of large $x$ limit and determinant.

Second we focus on the Fredholm determinant. The
key observation is that we can decompose the operator 
$|V|^{1/2}(\pa-A_0)^{-1}U|V|^{1/2}$
as follows. Direct computation using
the kernel $|V(x)|^{1/2}G(x;y)U(y)|V(y)|^{1/2}$, 
where $G$ is the kernel corresponding to $(\pa-A_0)^{-1}$, 
reveals  
\begin{equation*}
|V|^{1/2}(\pa-A_0)^{-1}U|V|^{1/2}=J-R,
\end{equation*}
where $J$ is the Volterra integral operator given in the
Volterra integral equation Lemma \ref{lemma:VIE} and $R$
is the integral operator with kernel given by 
\begin{equation*}
\phi_0^-(x)\hat Z_0^+(y)U(y)|V(y)|^{1/2}.
\end{equation*}
Note $R$ has separable kernel and is thus a finite rank operator
and trace class. Since we have 
$\|J\|_{\fJ_1}\leqslant
\bigl\||V|^{1/2}(\pa-A_0)^{-1}U|V|^{1/2}\bigr\|_{\fJ_1}+\|R\|_{\fJ_1}$
and $|V|^{1/2}(\pa-A_0)^{-1}U|V|^{1/2}$ 
is trace class by assumption, we deduce $J$ is also trace class. 
Using the product decomposition 
$\id-J+R=(\id-J)\bigl(\id+(\id-J)^{-1}R$ and
provided $J$ and $(\id-J)^{-1}R$ are trace class,
we have 
\begin{align*}
\mathrm{det}_1(\id-J+R)
&=\mathrm{det}_1\Bigl((\id-J)\bigl(\id+(\id-J)^{-1}R\bigr)\Bigr)\\
&=\mathrm{det}_1(\id-J)\,\mathrm{det}_1\bigl(\id+(\id-J)^{-1}R\bigr).
\end{align*}
As we see presently $(\id-J)^{-1}R$ is trace class 
as it has a separable kernel.

Third we compute $\mathrm{det}_1(\id-J)$. 
Using the relation $\mathrm{det}\,\exp=\exp\,\mathrm{tr}$ we have 
\begin{equation*}
\log\mathrm{det}_1(\id-J)=
-\sum_{\ell\geqslant1}\frac{1}{\ell}\,\mathrm{tr}\,J^\ell.
\end{equation*}
Since $J$ is a Volterra integral operator involving the Heaviside
function, using the trace power formula
Lemma~\ref{lemma:traceproductformula}, we observe that $\mathrm{tr}\,J^\ell$
will involve an integral over $\R^\ell$ of an integrand with
a factor $H(y_1-y_2)\,H(y_2-y_3)\cdots H(y_{\ell-1}-y_\ell)\,H(y_\ell-1)$
of a product of Heaviside functions.
Hence for all $\ell\geqslant 2$ the quantity $\mathrm{tr}\,J^\ell$ 
is zero as it is an integral of an integrand which is zero everywhere except 
on a subset of $\ell$-dimensional measure zero. We also observe 
$\mathrm{tr}\,J=0$. This follows using the unitary relation 
in Lemma~\ref{lemma:orthonormality} and the assumption 
$\mathrm{tr}\,V\equiv 0$ which imply the matrix trace of the kernel
of $J$ is continuous on the diagonal.
Hence we deduce $\mathrm{det}_1(\id-J)=1$.

Fourth we establish that $(\id-J)^{-1}R$ is trace class and compute 
$\mathrm{det}_1\bigl(\id+(\id-J)^{-1}R\bigr)$. 
Using that $\phi_0^-=(\mathrm{id}-J)\phi^-$ from the Volterra integral
equation Lemma~\ref{lemma:VIE}, we see 
\begin{align*}
\bigl((\id-J)^{-1}(R\varphi)\bigr)(x)
&=(\id-J)^{-1}\phi_0^-(x)\int_\R\hat Z_0^+(y)U(y)|V(y)|^{1/2}\varphi(y)\,\mathrm{d}y\\
&=\phi^-(x)\int_\R\hat Z_0^+(y)U(y)|V(y)|^{1/2}\varphi(y)\,\mathrm{d}y.
\end{align*}
Hence $(\id-J)^{-1}R$ has 
separable kernel $\phi^-(x)\hat Z_0^+(y)U(y)|V(y)|^{1/2}$ and is thus trace class.
Then using the result for separable kernels in Corollary~\ref{corollary:separablekernel},
we find for any $\epsilon>0$ sufficiently small,
\begin{align*}
\log\mathrm{det}_1\bigl(\id+\epsilon(\id-J)^{-1}R\bigr)
&=\sum_{\ell\geqslant1}\frac{(-1)^{\ell-1}}{\ell}\epsilon^\ell
\,\mathrm{tr}\,\bigl((\id-J)^{-1}R\bigr)^\ell\\
&=\sum_{\ell\geqslant1}\frac{(-1)^{\ell-1}}{\ell}\epsilon^\ell
\,\mathrm{tr}\,\biggl(\int_\R\hat Z_0^+(y)U(y)|V(y)|^{1/2}\phi^-(y)\,\mathrm{d}y\biggr)^\ell\\
&=\log\mathrm{det}\biggl(\id+\epsilon\int_\R (\hat Z_0^+VY^-)(y)\,\mathrm{d}y\biggr).
\end{align*}
By analytic continuation we can extend this result to $\epsilon=1$.
Then removing the logarithms, this equals our expression above for the 
transmission coefficient.
\qed
\end{proof}
\begin{remark}
We have the following observations on the results above:
(i) At the core of the equivalence theorem is, for solutions
decaying in the far field, on one hand we have the Fredholm determinant 
associated with the solvability of $(\id-J+R)\phi=O$, while on the other
we have the transmission coefficient associated with
the solvability of $(\id-J)\phi^-=\phi_0^-$;
(ii) We note we have $(\pa-A_0)|V|^{-1/2}J|V|^{1/2}=V$ and 
$(\pa-A_0)|V|^{-1/2}R|V|^{1/2}=O$;
(iii) In essence, in the final calculation in the proof
of Theorem~\ref{theorem:equivalence}, we demonstrate that the 
trace of any power of $(\id-J)^{-1}R$ equals the trace of the
corresponding power of $\int(\hat Z_0^+VY^-)(y)\,\mathrm{d}y$.
Hence we could prove the equivalence of the determinants 
using the Plemelj--Smithies formula for 
$\mathrm{det}_1\bigl(\id+\epsilon(\id-J)^{-1}R\bigr)$ which
is analytic in $\C$---see Simon~\cite[Theorem~6.8]{BS77};
(iv) Since we can swap the order of two arguments under the trace, using
the trace formula Lemma~\ref{lemma:traceproductformula}, 
we observe the traces of all powers of $|V|^{1/2}(\pa-A_0)^{-1}U|V|^{1/2}$ 
and $(\pa-A_0)^{-1}V$ coincide. Hence the Fredholm determinants
of the two operators coincide, assuming the Fredholm determinant
of $\id-(\pa-A_0)^{-1}V$ exists; 
(v) If $\mathrm{tr}\,V\neq0$, then we need to include 
the factor $\exp(-\mathrm{tr}\,J)$ in the evaluation of the 
Fredholm determinant;
(vi) The approach we used in the proof of the equivalence
theorem is based on the standard approach---decomposing 
the given operator with semi-separable kernel into the sum
of a Volterra and finite rank operators---which is given  
in Gesztesy, Latushkin \& Makarov~\cite{GL07} for Hilbert--Schmidt 
operators. They utilize results from Gesztesy \& Makarov~\cite{GK04},
Gohberg, Goldberg \& Kaashook~\cite{ISM90} and
Gohberg, Goldberg \& Krupnik~\cite{GGK}; and 
(vii) Gesztesy \textit{et al.\/ }~\cite{GL07} do not assume strict 
hyperbolicity which we have done to keep the arguments as 
succinct as possible.
\end{remark}

\section{\textsf{Distinct far fields}}\label{sec:distinctfarfields}
Our goal in this section is to show that the equivalence of the
transmission coefficient and Fredholm determinant carries over
to the case when the far field limits of the eigenvalue problem
$(\pa-A_0-V)Y=O$ are distinct. The Evans function and transmission
coefficient are well defined in this instance with only minor
modification to their construction in Section~\ref{sec:dets}---which 
is the standard approach. Indeed what underlies this approach is
that we decompose $A_0+V$ in such a way that $A_0=A_0(x)$ ensures
the distinct far field limits are satisfied so $V\to O$ as $x\to\pm\infty$.
We assume this decomposition in this section, as it
also gives us the natural framework to show the equivalence 
of the transmission coefficient and a Fredholm determinant with only a slight 
modification to the arguments in the last section. 
We also assume in each far field $A_0$ is constant and strictly hyperbolic, 
with the same hyperbolic splitting (the number of eigenvalues with 
positive real parts, say $k$, and negative real parts, consequently $n-k$).
We can thus also assume $V\to O$ as $x\to\pm\infty$.
Further, we assume both $A_0$ and $V$ are continuous and uniformly bounded. 
To construct the Evans function or transmission coefficient we do not
need to perform this decomposition, however the properties
we derive for the solutions of the equations
\begin{equation*}
\pa Y_0=A_0Y_0
\qquad\text{and}\qquad
\pa Z_0=-Z_0A_0, 
\end{equation*}
are crucial to our equivalence proof. We now have two separate 
constant coefficient equations in the far field 
corresponding to $\pa Y_0=A_0Y_0$. However as before in 
Section~\ref{sec:dets}, given the identical hyperbolic splitting
of the two far field limits of $A_0$, there is $k$-dimensional subspace of 
solutions that decays exponentially as $x\to-\infty$, and an 
$(n-k)$-dimensional subspace of solutions that decays exponentially as $x\to+\infty$.
We suppose these subspaces are given by the column span of solutions
$Y_0^-\in\C^{n\times k}$ and $Y_0^+\in\C^{n\times(n-k)}$, respectively.
For the adjoint equation $\pa Z_0=-Z_0A_0$ there exist
solutions collected as rows in the matrices $Z_0^-\in\C^{(n-k)\times n}$ 
and $Z_0^+\in\C^{k\times n}$ which decay as $x\to-\infty$ and $x\to+\infty$, 
respectively. Solutions of commensurate dimension $Y^\pm$ exist to the full 
problem $(\pa-A_0-V)Y=O$ such that $Y^\pm\sim Y_0^\pm$ as $x\to\pm\infty$.

We can recover the unitary relations  
of Lemma~\ref{lemma:orthonormality} for $Y_0^\pm$ and 
suitably defined scaled solutions $\hat Z_0^\pm$
despite distinct far fields. The diagonal relation
for $Y_0^\pm$ and $Z_0^\pm$ becomes the following.
\begin{lemma}[\textbf{\textsf{Diagonal block relation}}]\label{lemma:DBR}
The solutions $Y_0^\pm$ and $Z_0^\pm$ satisfy the relation
\begin{equation*}
\begin{pmatrix}
Z_0^+\\ Z_0^-
\end{pmatrix}
\begin{pmatrix}
Y_0^- & Y_0^+
\end{pmatrix}=\begin{pmatrix}D_-&O\\O&D_+\end{pmatrix},
\end{equation*}
where $D_\pm$ are constant matrices. 
\end{lemma}
\begin{proof}
We observe, by direct computation for any pair of 
solutions $Y_0\in\C^{n\times 1}$ and $Z_0\in\C^{1\times n}$, that 
$\pa (Z_0Y_0)=O$ and so $Z_0Y_0$ is constant. Asymptotically 
we have $Y_0^-\sim\exp\bigl(A_0(-\infty)x\bigr)U^-$ and 
$Z_0^-\sim W^-\exp\bigl(-A_0(-\infty)x\bigr)$ where the columns of $U^-$
and rows of $W^-$ are the right and left eigenvectors of $A_0(-\infty)$.
Similarly we have $Y_0^+\sim\exp\bigl(A_0(+\infty)x\bigr)U^+$ and 
$Z_0^+\sim W^+\exp\bigl(-A_0(+\infty)x\bigr)$ where the columns of $U^+$
and rows of $W^+$ are the right and left eigenvectors of $A_0(+\infty)$.
Then using that $W^-U^-=O$ and $W^+U^+=O$, we conclude 
$Z_0^-Y_0^-=O$ and $Z_0^+Y_0^+=O$. The remaining constant matrices
are thus $Z_0^+Y_0^-=D_-$ and $Z_0^-Y_0^+=D_+$.
\qed
\end{proof}
Note that the constant matrices $D_\pm$ are in general no longer diagonal 
as $Y_0^-$ and $Z_0^-$, and, $Y_0^+$ and $Z_0^+$, satisfy different 
asymptotic problems. We must keep in mind that $Y_0^\pm$ and $Z_0^\pm$ 
are the solutions generated by $A_0=A_0(x)$. We \emph{assume} the 
constant matrices $Z_0^+Y_0^-=D_-$ and $Z_0^-Y_0^+=D_+$ are \emph{nonsingular}.
Generically this will be the case---see the discussion in the conclusions
Section~\ref{sec:conclu}. This guarantees the free transmission coefficient
is non-zero. Using the diagonal block relation Lemma~\ref{lemma:DBR},
it guarantees the free Evans function is non-zero. This guarantees 
$(\pa-A_0)^{-1}$ exists (this was automatic in the equal far field case
under the strict hyperbolicity assumption on the constant matrix $A_0$). 
With this proviso, the block diagonal relation implies 
the framework we provided in Section~\ref{sec:dets} for establishing  
the relation between the Evans function and the transmission
coefficient carries through essentially unchanged. We can mirror
all the results in Section~\ref{sec:dets} up to and 
including the schematic formula between the Evans function 
and transmission coefficient.
Then again with the same proviso, the scaled unitary solutions 
are defined as for the equal far field case in Section~\ref{sec:dets} by
$\hat W^\pm\coloneqq D_\mp^{-1}W^\pm$ and $\hat Z_0^\pm\coloneqq D_\mp^{-1}Z_0^\pm$.
\begin{lemma}[\textbf{\textsf{Unitary relations reprise}}] 
Unitary relations for $Y_0^\pm$ and $\hat Z_0^\pm$ hold for distinct far fields.
\end{lemma}
\begin{proof}
Making the corresponding change of variables in the diagonal
block relation above generates the same diagonal block relation
but with $\hat Z_0^\pm$ in place of $Z_0^\pm$ and the appropriate 
identity matrices in place of $D_\pm$ on the right. The two 
matrices on the left of the new diagonal block relation are then 
inverses of each other and the unitary relation in reverse order follows.
\qed
\end{proof}
The framework we provided in Section~\ref{sec:equivtheorem} for 
establishing the relation between the transmission
coefficient and the Fredholm determinant 
$\mathrm{det}_1\bigl(\id-|V|^{1/2}(\pa-A_0)^{-1}U|V|^{1/2}\bigr)$, where now $A_0=A_0(x)$, 
carries through with only a couple of modifications which we now outline. 
The Green's function $G=G(x,y)$ associated with $(\pa-A_0)^{-1}$ has
the same semi-separable form except that $Y_0^\pm$ and $\hat Z_0^\pm$ 
are now the solutions generated by $A_0=A_0(x)$. Again the 
unitary relations imply the jump condition for $G$ across the diagonal line $y=x$
is satisfied. If we assume $\mathrm{tr}\,V\equiv 0$ then 
the trace of $|V|^{1/2}(\pa-A_0)^{-1}U|V|^{1/2}$ or $(\pa-A_0)^{-1}V$ 
is defined in the same way.
The Volterra integral equation $\phi^-=\phi_0^-+J\phi^-$ of Lemma~\ref{lemma:VIE}
still applies but now only for $J$ defined via the kernel 
$H(x-y)\bigl(\phi_0^-(x)Z_0^+(y)+\phi_0^+(x)Z_0^-(y)\bigr)U(y)|V(y)|^{1/2}$. 
This requires independent proof which does not rely on $A_0$ being
constant as follows. Direct computation using the unitary relations and that
$(\pa-A_0)Y_0^\pm=O$, reveals  
\begin{align*}
(\pa-A_0)\bigl(Y_0^-+|V|^{-1/2}&J|V|^{1/2}Y^-\bigr)(x)\\
&=(\pa-A_0)\bigl(|V|^{-1/2}J|V|^{1/2}Y^-\bigr)(x)\\
&=(\pa-A_0)\int_{-\infty}^x\bigl(Y_0^-(x)\hat Z_0^+(y)+Y_0^+(x)
\hat Z_0^-(y)\bigr)V(y)Y^-(y)\,\mathrm{d}y\\
&=\bigl(Y_0^-(x)\hat Z_0^+(x)+Y_0^+(x)\hat Z_0^-(x)\bigr)V(x)Y^-(x)\\
&=V(x)Y^-(x).
\end{align*}
Hence if $Y^-=Y_0^-+|V|^{-1/2}J|V|^{1/2}Y^-$ 
we have just shown that $(\pa-A_0)Y^-=VY^-$.
Premulti-plying this Volterra equation for $Y^-$ by $|V|^{1/2}$
and using the definitions for $\phi^-$ and $\phi^-_0$ gives the result. 
The proof of the equivalence Theorem~\ref{theorem:equivalence} now follows 
step by step with the solutions $Y_0^\pm$ and $\hat Z_0^\pm$ 
those generated by $A_0=A_0(x)$. 
We summarize these conclusions as follows.
\begin{theorem}[\textbf{\textsf{Equivalence reprised}}]
Assume 
$\pa-A_0$ is invertible, 
$|V|^{1/2}(\pa-A_0)^{-1}U|V|^{1/2}$ is trace class,
$\mathrm{tr}\,V\equiv 0$, $A_0=A_0(x)$ and $V=V(x)$ are both continuous
and uniformly bounded and $V=V(x)$ is $w^2$-integrable. 
Further assume in each far field $A_0$ is constant 
and strictly hyperbolic, with the same hyperbolic splitting, and 
$V\to O$ as $x\to\pm\infty$. Then the transmission coefficient 
for the eigenvalue problem $(\pa-A_0-V)Y=O$ and Fredholm determinant 
of $\mathrm{id}-|V|^{1/2}(\pa-A_0)^{-1}U|V|^{1/2}$ are equal.
\end{theorem}

\section{\textsf{Conclusions}}\label{sec:conclu}
We have established the equivalence of the Evans function and transmission
coefficient for the eigenvalue problem $(\pa-A_0-V)Y=O$, 
in the sense that the ratio of the Evans function to
the free Evans function is equal to the ratio of the transmission 
coefficient to the free transmission coefficient. As we remarked
at the end of Section~\ref{sec:dets}, the Evans function and free
Evans function are invariant to the unitary rescaling, as is
the ratio of the transmission coefficient to the free transmission coefficient.
Hence the Fredholm determinant 
$\mathrm{det}_1\bigl(\id-|V|^{1/2}(\pa-A_0)^{-1}U|V|^{1/2}\bigr)$,
which equals the transmission coefficient with unitarily scaled solutions
for which the free transmission coefficient is unity, equals the ratio
of the transmission coefficient to the free transmission coefficient---whether
the unitary scaling is used or not. These statements hold for 
equal or distinct far fields. In other words we have established that
\begin{equation*}
\frac{\text{Evans function}}{\text{Free Evans function}}
=\frac{\text{Transmission coefficient}}{\text{Free transmission coefficient}}
=\text{Fredholm determinant}.
\end{equation*}
Let us now pull back our perspective to our original
problem of determining values of $\lambda\in\C$ for which there exist
solutions to $(\pa-A_0-V)Y=O$ where $A_0=A_0(x;\lambda)$ in general. 
The locale of the essential spectrum is determined by the values of $\lambda$
for which the far field limits of $A_0$ are no longer strictly hyperbolic---at 
least one eigenvalue of either limit becomes pure imaginary characterizing the 
continuous spectrum or the hyperbolic splitting no longer matches. 
Away from the essential spectrum, we must choose $A_0=A_0(x;\lambda)$ 
suitably so that $(\pa-A_0)^{-1}$ exists---as mentioned previously this
is not an issue in the equal far field case. Hence, away from the essential
spectrum, the free Evans function and free transmission coefficients are 
bounded and non-zero. Since the Evans function is analytic in that region,
the product of the free Evans function and the ratio of the transmission
and free transmission coefficients is analytic there, as well as the product
of the free Evans function and the Fredholm determinant. Zeros of these
product quantities must thus coincide in this region. Those zeros coincide
with pure point eigenvalues with coincident multiplicity.  
We remark that for eigenvalue problems of the 
form $(\pa-A_0-V)Y=O$ that arise in the study of the stability of travelling
waves the origin is an eigenvalue associated with translation invariance.
In some cases the origin is embedded in the essential spectrum within 
which further analysis of all the discriminants above is required. 
Some instructive explicit examples
illustrating the relation, between the Evans function and transmission
coefficient can be found in Kapitula \& Sandstede~\cite{KS04} and Kapitula~\cite{K}, 
between the transmission coefficient and Fredholm determinant can be found in
Simon~\cite[p.~51]{BS05} and between the Evans function and $2$-modified Fredholm
determinant in Gesztesy \textit{et al.\/ } \cite{GL07} 
and Gesztesy \textit{et al.\/} \cite{GYZ08}. 

For the eigenvalue problem $(\pa-A_0-V)Y=O$ where $A_0=A_0(x;\lambda)$ and 
$\lambda\in\C$ is the spectral parameter, our interest surrounds families
of operators $\pa-A_0-V$ parameterized by $\lambda$. For example,
the analyticity of the Evans function means we can conduct a global search 
for eigenvalues via contour integration and the residue theorem in any subregion
away from the essential spectrum. Here the family of operators 
consists of those associated with values of $\lambda$ parameterizing 
the boundary contour of the subregion. Hence the 
object of interest is a determinant line bundle. See Quillen~\cite{Q}, 
Jost~\cite[Section~6.9]{Jost}, Deng \& Nii~\cite{DN} and Alexander, 
Gardner and Jones~\cite{AGJ90} for constructions of such line bundles. 
In the Alexander \emph{et al.} line bundle construction, the fibres are 
explicitly the Evans function. These constructions suggest we 
define, away from the essential spectrum, the determinant of the unbounded Fredholm
operator $\pa-A_0-V$ by the Evans function corresponding to the 
problem $(\pa-A_0-V)Y=O$, i.e.\/ we set
\begin{equation*}
\mathrm{det}(\pa-A_0-V)\coloneqq\text{Evans function},
\end{equation*}
or for a coordinate free definition, as the ratio of the Evans function
to the free Evans function. 
With the definition above, the determinant $\det(\pa-A_0)$ is the free Evans
function. And in this paper, assuming $\mathrm{tr}V\equiv0$, 
we have established that
\begin{equation*}
\mathrm{det}(\pa-A_0-V)
=\mathrm{det}(\pa-A_0)
\cdot\mathrm{det}_1\bigl(\id-(\pa-A_0)^{-1}V\bigr),
\end{equation*}
representing the natural determinant multiplicative property for such operators.
Elliptic operators for which this does not hold, i.e.\/ when there are 
multiplicative anomalies, is a research field in itself with important 
applications in mathematical physics. Here the zeta-regularized,
Quillen, Segal and regularized Fredholm determinants and their relations
are studied. See Quillen~\cite{Q}, Kontsevich \& Vishik~\cite{KV94} 
and Scott \& Wojciechowski~\cite{SW} for more details.

Finally, in the study of the linear stability of 
travelling wave solutions to nonlinear partial differential equations
often numerical simulation is required in the search for eigenvalues. 
Indeed one of the main motivations for establishing the 
connection between the Evans function and Fredholm determinant was the need
to compute the stability of multidimensional travelling waves; see
for example Sandstede \& Scheel~\cite{BA08}. 
Gesztesy, Latushkin \& Zumbrun~\cite[Section~4]{GYZ08} 
established important convergence results for suitable Galerkin 
approximations for the multidimensional problem. Humpherys \&
Zumbrun~\cite{HZ06}, Ledoux, Malham \& Th\"ummler~\cite{LMT10} and
Ledoux, Malham, Niesen \& Th\"ummler~\cite{LMNT10} have also designed
numerical methods for computing the Evans function for the multidimensional case. 
The computation of Fredholm determinants also naturally arise in the computation
of the Maslov index for linear symplectic systems, in particular in the 
multi-dimensional case which relies on computing flows in the 
Fredholm Lagrangian Grassmannian manifold, see Deng \& Jones~\cite{DJ}
and Beck \& Malham~\cite{BM}.
Finally, a natural question to ask that is still open is, as numerical
tools, which of the Evans function, transmission coefficient, Fredholm determinant 
or Galerkin-type direct projection methods has the optimal complexity?
See Karambal~\cite{IK13} for some results in this direction. 
A comprehensive study of this complexity issue would be a useful
follow-on project.

\section*{\textsf{Acknowledgement}}
We thank the referees for their very helpful hints, comments
and suggestions that have substantially improved our original manuscript. 
We also thank the third referee who lead us to consider the determinant of
the unbounded operator in our Conclusions also as 
the coordinate free ratio mentioned.
IK would also like to thank the Numerical Algorithms and Intelligent Software Centre 
funded by the UK EPSRC grant EP/G036136 and the Scottish Funding Council. 


\appendix


\section{\textsf{Proof of the Hilbert--Schmidt class lemma}}\label{app:HSlemma}

\begin{proof}
First, given $K\in\fJ_\infty$, 
the existence of a kernel function $G\in L^2(\R^2;\C^{n\times n})$ 
if and only if $K\in\fJ_2$, follows from standard existence theory; 
see Reed \& Simon~\cite[p.~210]{RSI}.
Second, we prove the isometry of the map $G\mapsto K$ from
$L^2(\R^2;\C^{n\times n})$ to $\fJ_2$. Suppose the set of $\C^n$-valued
functions $\{\varphi_m\}_{m\geqslant1}$ are an unitary basis for $L^2(\R;\C^n)$.
The set $\{\varphi_\ell\otimes\varphi_m^\dag\}_{\ell,m\geqslant1}$ is
a unitary basis for the Hilbert space $L^2(\R^2;\C^{n\times n})$.
Hence for any $G\in L^2(\R^2;\C^{n\times n})$, there exists a double sequence 
of constants $G_{\ell,m}\in\C^{n\times n}$ whose trace is square-summable, such that
$G(x;y)=\sum_{\ell,m\geqslant1}G_{\ell,m}\varphi_\ell(x)\varphi_m^\dag(y)$.
By direct calculation, we have
\begin{align*}
\mathrm{tr}\,|K|^2&=\sum_{m\geqslant1}\la\varphi_m,K^\dag K\varphi_m\ra_{L^2(\R;\C^n)}\\
&=\sum_{m\geqslant1}\la K\varphi_m,K\varphi_m\ra_{L^2(\R;\C^n)}\\
&=\sum_{m\geqslant1}\int_{\R}(K\varphi_m)^\dag(x)(K\varphi_m)(x)\,\mathrm{d}x\\
&=\sum_{m\geqslant1}\int_{\R^3}\bigl(G(x,\xi)\varphi_m(\xi)\bigr)^\dag
\bigl(G(x,\eta)\varphi_m(\eta)\bigr)\,\mathrm{d}\xi\,\mathrm{d}\eta\,\mathrm{d}x\\
&=\sum_{m\geqslant1}\int_{\R^3}\varphi_m^\dag(\xi)G^\dag(x,\xi)
G(x,\eta)\varphi_m(\eta)\,\mathrm{d}\xi\,\mathrm{d}\eta\,\mathrm{d}x.
\end{align*}
Using the expansion for $G=G(x;y)$ above 
we find, 
\begin{align*}
\mathrm{tr}\,|K|^2
&=\sum_{m,k,\ell,p,q\geqslant1}\int_{\R^3}\varphi_m^\dag(\xi)
\bigl(G_{k,\ell}\varphi_k(x)\varphi_\ell^\dag(\xi)\bigr)^\dag
\bigl(G_{p,q}\varphi_p(x)\varphi_q^\dag(\eta)\bigr)\varphi_m(\eta)
\,\mathrm{d}\xi\,\mathrm{d}\eta\,\mathrm{d}x\\
&=\sum_{m,k,\ell,p,q\geqslant1}\int_{\R^3}\varphi_m^\dag(\xi)\varphi_\ell(\xi)\varphi_k^\dag(x)
G_{k,\ell}^\dag G_{p,q}\varphi_p(x)\varphi_q^\dag(\eta)\varphi_m(\eta)
\,\mathrm{d}\xi\,\mathrm{d}\eta\,\mathrm{d}x\\
&=\sum_{m,k,p\geqslant1}\int_{\R}\varphi_k^\dag(x)G_{k,m}^\dag G_{p,m}\varphi_p(x)\,\mathrm{d}x\\
&=\mathrm{tr}\,\sum_{m,k,p,r\geqslant1}\int_{\R^2}\varphi_{m}(y)
\varphi_k^\dag(x)G_{k,m}^\dag G_{p,r}\varphi_p(x)\varphi_r^\dag(y)\,\mathrm{d}x\,\mathrm{d}y\\
&=\mathrm{tr}\,\sum_{m,k,p,r\geqslant1}\int_{\R^2}\bigl(G_{k,m}\varphi_k(x)\varphi_m^\dag(y)\bigr)^\dag
\bigl(G_{p,r}\varphi_p(x)\varphi_r^\dag(y)\bigr)\,\mathrm{d}x\,\mathrm{d}y\\
&=\int_{\R^2}\mathrm{tr}\,G^\dag(x;y)G(x;y)\,\mathrm{d}x\,\mathrm{d}y\\
&=\|G\|_{L^2(\R^2;\C^{n\times n})}^2.
\end{align*}
\qed
\end{proof}


\section{\textsf{Proof of the trace formula lemma}}\label{app:traceformula}

\begin{proof}
For any $\ell\in\mathbb N$, suppose that $K_1$, $K_2$, \ldots, $K_\ell$ 
are Hilbert--Schmidt operators with canonical respective kernels $G_1$, $G_2$, \ldots, $G_\ell$.
As in the proof of the Hilbert--Schmidt class lemma 
above, there exists a double sequence of constants $G_{p,m}^{(\ell)}\in\C^{n\times n}$ 
whose trace is square-summable, such that
$G_\ell(x;y)=\sum_{p,m\geqslant1}G_{p,m}^{(\ell)}\varphi_p(x)\varphi_m^\dag(y)$.
Then by direct computation we see that  
\begin{align*}                                                                                    
\mathrm{tr}&\,K_1\cdots K_\ell  \\    
&=\sum_{m\geqslant1}\la\varphi_m,K_1\cdots K_\ell\varphi_m\ra_{L^2(\R;\C^n)}\\  
&=\sum_{m\geqslant1}\int_{\R}\varphi_m^\dag(x)\bigl(K_1\cdots K_\ell
\varphi_m\bigr)(x)\,\mathrm{d}x\\                                                                 
&=\sum_{m\geqslant1}\int_{\R^{\ell+1}}\varphi_m^\dag(x)                                           
G_1(x;y_1)G_2(y_1;y_2)\cdots G_\ell(y_{\ell-1};y_\ell)\varphi_m(y_\ell)         
\,\mathrm{d}y_\ell\cdots\,\mathrm{d}y_1\,\mathrm{d}x\\                                            
&=\sum_{m,p,q\geqslant1}\int_{\R^{\ell+1}}\varphi_m^\dag(x)                                   
G_1(x;y_1)\cdots \\
&\qquad\qquad\qquad\qquad\cdots G_{\ell-1}(y_{\ell-2};y_{\ell-1})
G_{p,q}^{(\ell)}\varphi_p(y_{\ell-1})\varphi^\dag_q(y_\ell)\varphi_m(y_\ell)            
\,\mathrm{d}y_\ell\cdots\,\mathrm{d}y_1\,\mathrm{d}x\\                                            
&=\sum_{m,p\geqslant1}\int_{\R^{\ell}}\varphi_m^\dag(x)                                           
G_1(x;y_1)\cdots G_{\ell-1}(y_{\ell-2};y_{\ell-1})                        
G_{p,m}^{(\ell)}\varphi_p(y_{\ell-1})                                     
\,\mathrm{d}y_{\ell-1}\cdots\,\mathrm{d}y_1\,\mathrm{d}x\\                                        
&=\mathrm{tr}\,\int_{\R^{\ell}}\sum_{m,p\geqslant1}                                               
G_1(x;y_1)\cdots G_{\ell-1}(y_{\ell-2};y_{\ell-1})             
G_{p,m}^{(\ell)}\varphi_p(y_{\ell-1})\varphi_m^\dag(x)                 
\,\mathrm{d}y_{\ell-1}\cdots\,\mathrm{d}y_1\,\mathrm{d}x\\                                        
&=\mathrm{tr}\,\int_{\R^{\ell}}                                                                   
G_1(x;y_1)G_2(y_1;y_2)\cdots G_{\ell-1}(y_{\ell-2};y_{\ell-1})              
G_\ell(y_{\ell-1};x)\,\mathrm{d}y_{\ell-1}\cdots\,\mathrm{d}y_1\,\mathrm{d}x\\                         
&=\mathrm{tr}\,\int_{\R^{\ell}}                                                                   
G_1(y_1;y_2)G_2(y_2;y_3)\cdots G_{\ell-1}(y_{\ell-1};y_{\ell}) 
G_\ell(y_{\ell};y_1)\,\mathrm{d}y_{\ell}\cdots\,\mathrm{d}y_1.       
\end{align*}   
The proof for the case of a single trace class operator $K$ with continuous
kernel $G$ exactly follows the argument above with $\ell=1$.
\qed
\end{proof}

\end{document}